\documentclass[final]{siamltex}
\usepackage{amsmath}
\usepackage{amsfonts}
\usepackage{amssymb}
\usepackage{graphicx}
\usepackage{url}
\usepackage{color}
\usepackage{fullpage}
\usepackage{subfigure,float}
\usepackage{bm}
\usepackage{algorithmic}
\usepackage{algorithm}
\usepackage{color}
\usepackage{epstopdf}
\begin{document}
\baselineskip 0.58cm 
\title{Nonlinear Filters for Hidden Markov Models of Regime Change with Fast Mean-Reverting States\footnote{This research has been partially supported by NSF grant \#DMS-0739195.}}
\author{Andrew Papanicolaou\footnote{Department of Operations Research and Financial Engineering, Princeton University, Princeton NJ 08544, apapanic@princeton.edu}}
\maketitle
\begin{abstract}We consider filtering for a hidden Markov model that evolves with multiple time scales in the hidden states. In particular, we consider the case where one of the states is a scaled Ornstein-Uhlenbeck process with fast reversion to a shifting-mean that is controlled by a continuous time Markov chain modeling regime change. We show that the nonlinear filter for such a process can be approximated by an averaged filter that asymptotically coincides with the true nonlinear filter of the regime-changing Markov chain as the rate of mean reversion approaches infinity. The asymptotics exploit weak converge of the state variables to an invariant distribution, which is significantly different from the strong convergence used  to obtain asymptotic results in \cite{papanicolaou2010}.
\end{abstract}

\section{Introduction}
\label{sec:into}
When considering systems of stochastic differential equations with multiple time scales, hidden Markov models (HMMs) come up because we often want to include regime changes. In estimating the hidden states, filtering methods are used to obtain posterior distributions. If fast mean-reversion is involved, one might expect that the posterior distribution of the mean-reverting states is tempered only by an invariant measure and not by the history of the data. This is indeed correct, provided that observations are sampled discretely, allowing for relaxation in the distribution of the fast mean-reverting states by the time new observations arrive. Thus, it is reasonable to approximate the filtering distribution with an asymptotic average over these states.

An asymptotic nonlinear filter can be useful because it may be of reduced dimension, and therefore considerably easier to compute. By taking a Bayesian approach to filtering, the posterior distributions are computed by integrating exhaustively over all states. These integrations are prohibitively slow for systems with hidden states of multiple dimension. In contrast, the computational demands for filters of linear systems with Gaussian white noise are usually not an issue because we use Kalman filters. They are obtained from projections onto the linear subspace of observable data, while filtering in the presence of nonlinearity and/or non-Gaussian components does not allow this linear theory. However, it is possible to find faster nonlinear filtering algorithms when the HMM has states whose law quickly converges to an invariant measure. Such filters are faster to compute because it becomes unnecessary to track the hidden states where convergence quickly occurs.

Asymptotic analysis of multi-scale stochastic models has been done extensively by Yin and Zhang in \cite{yinZhangCont,yinZhangDiscrete}, in which they also analyze the asymptotics of the Kalman filter and the Wonham filter. Multi-scale asymptotics are also of interest in stochastic volatility models \cite{FPS00,howison}, which suggests that filtering is also useful in finance. In the setting of multiple time-scales (in particular HMMs with fast mean-reverting states), nonlinear filtering has, to our knowledge, not been explored to its full extent. The basic nonlinear filtering theory is presented by Rabiner \cite{rabiner1989}, Jazwinski \cite{jazwinski1970}, and Yin and Zhang \cite{yinZhangDiscrete} as a fundamental Bayesian method to track the hidden states of nonlinear HMMs. Nonlinear filtering has been used extensively in target-tracking \cite{barShalom,petrov2000}, but in these types of problems it may be more efficient to generate particles for the nonlinear components of the state space and then exploit the Kalman filter everywhere else, as is done in \cite{gustafsson2002,gustafsson2005} where they propose the use of marginalized filters which achieve a Rao-Blackwell lower bound given the sufficient statistics provided by the Kalman filter. Essentially, nonlinear filtering of continuous state-space random processes needs to be numerically computed with some kind of approximating scheme. An important result in this area is that of Kushner \cite{kushner08}, which shows consistency and stability for Markov chain approximations of the nonlinear filter with continuous time observations. In general, particle filters \cite{asmussen2007,cappe2005} are a method used frequently for approximating the nonlinear filter, but the amount of time required to generate the necessary Monte Carlo samples can be a deterrent in practice. 

In this paper, we consider in detail a nonlinear multi-scale model with a mean-reverting state that reverts faster as a parameter $\epsilon>0$ gets smaller. The model has applications in target-tracking in finance and elsewhere. In the early sections we discuss how a distributional averaging occurs in the posterior as $\epsilon\searrow0$ if the data is sampled discretely, and we also give some examples to illustrate some possible variations of the model. We then switch to a rigorous analysis and provide a detailed proof showing that indeed an averaged filter is obtained as $\epsilon\searrow 0$, provided that appropriate assumptions are made when defining the model. We also present some numerical experiments to illustrate how such average filters perform compared to other approximations of the nonlinear filter. 

While the theory presented is for a specific model, the proofs can be carried out for other, more general process of a similar nature. \textcolor{black}{This paper will focus on models where the fast mean-reverting state is a one-dimensional Ornstein-Uhlenbeck process, but the theory can be reworked for other diffusion such as Cox-Ingersol-Ross processes. Also, the theory is easily generalized to the multivariate setting, both for cases when observations are a vector and when the OU is a vector. In fact, a model with a multivariate OU process and a scalar-valued regime process can be well-managed by our filter of reduced-dimensionality if the asymptotic  approximations apply. Parameter estimation is also an issue of interest, but one which we do not cover in this paper. The expectation maximization (EM) algorithm is implementable for the asymptotic model and will use the asymptotic filter that is computed in this paper, however the asymptotics of EM are perhaps a topic for future research.}

This paper is organized as follows: in Section \ref{sec:formulation} we present the HMM which is considered in this paper; the associated nonlinear filter is written in Section \ref{sec:filter}; the asymptotically averaged filter is given in Section \ref{sec:heuristics}; some example applications are given in Section \ref{sec:examples}; and Section \ref{sec:issues} describes some particular effects that scaling and observation sampling can have on the asymptotic filter. In Section \ref{sec:rigorousProof} we restate the model and the problem in a mathematically rigorous manner, and then goes on to a proof, while relegating to Appendix \ref{sec:details} some of the details that are included for the sake of completeness. Finally, in Section \ref{sec:numerics} we give a detailed account of a numerical simulation to illustrate how the average filter performs relative to some of the parameter choices.

\section{Formulation}
\label{sec:formulation}
In this section we introduce the model that is the focus of the paper. The notation of sections \ref{sec:process} and \ref{sec:filter} is used throughout the paper. In Section \ref{sec:heuristics} we discuss the averaged filter that occurs as the rate of mean reversion approaches infinity, in \ref{sec:examples} we give some examples, and in \ref{sec:issues} we consider changes in the averaged filter when using a somewhat different scaling. 
\subsection{The Process}
\label{sec:process}
\begin{figure}[htbp] 
   \centering
   \includegraphics[width=5.5in]{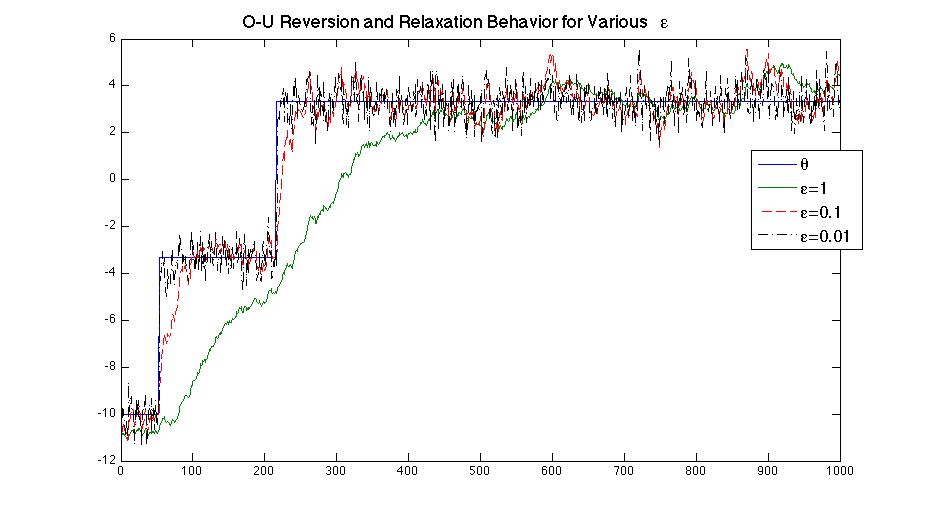} 
   \caption{\small An example of the proposed OU-type process with $Q(x)\equiv Q$ (i.e. $Q$ is a constant indendent of $x$). In this example we see the stationary behavior of $X^\epsilon(t)$ around $\Theta(t)$ occur faster and faster as $\epsilon$ decreases.}
   \label{fig:OUexample}
\end{figure}

Given the parameter $\epsilon>0$, which denotes the mean reversion time, we let $t\geq0$ denote time and we consider a hidden Markov model (HMM) consisting of the following processes:
\begin{eqnarray}
\nonumber
\Theta^\epsilon(t)\qquad&=&\qquad\hbox{a hidden Markov chain in a finite state space}\\
\nonumber
X^\epsilon(t)\qquad&=&\qquad\hbox{a hidden diffusion with mean reverting properties}\\
\nonumber
Y^\epsilon(t)\qquad&=&\qquad\hbox{an observed diffusion process.}
\end{eqnarray}
The hidden pair $(X^\epsilon(t),\Theta^\epsilon(t))\in\mathbb R\times \mathcal S$ is a Markov process with $X^\epsilon(t)$ being a real-valued process and $\Theta^\epsilon(t)$ taking value in the finite discrete space $\mathcal S\doteq\{s_1,\dots,s_M\}$. The forward Kolmogorov equation for $\Theta^\epsilon(t)$ is written in terms of a transition intensity matrix $Q\in\mathbb R^{M\times M}$ that is a function of $X^\epsilon(t)$,
\begin{eqnarray}
\label{eq:dPtheta}
\frac{\partial}{\partial t}\mathbb P(\Theta^\epsilon(t)=s_i)&=&\sum_jQ_{ji}(X^\epsilon(t))\mathbb P(\Theta^\epsilon(t)=s_j)~,~~ t>0\\
\mathbb P(\Theta^\epsilon(0)=s_i)&=&\rho_0(s_i)
\end{eqnarray}
for all $i\leq M$. We assume that jumps in $\Theta^\epsilon(t)$ cannot occur at an infinite rate, and that there are no cemetery states. Therefore there are constants $\beta$ and $\alpha$ such that
\begin{eqnarray}
\label{eq:upperJumpRateBound}
\sup_x(-Q_{ii}(x))&\leq& \beta<\infty\\
\label{eq:lowerJumpRateBound}
\inf_x(-Q_{ii}(x))&\geq&\alpha>0
\end{eqnarray}
for all $i\leq M$ and all $x\in \mathbb R$. The process $X^\epsilon(t)$ is a type of Ornstein-Uhlenbeck (OU) with $\Theta^\epsilon(t)$ serving as a regime changing shifting mean,
\begin{equation}
\label{eq:dX}
dX^\epsilon(t) = \frac{1}{\epsilon}(\Theta^\epsilon(t)-X^\epsilon(t))dt+\frac{1}{\sqrt\epsilon}dW(t)
\end{equation}
where $W(t)$ is an independent Wiener process and the law of $X^\epsilon(0)$ is given. We say that `$1/\epsilon$' is the rate of mean-reversion for $X^\epsilon(t)$. This model is a regime-switching OU, and as $\epsilon\searrow 0$ there is fast relaxation to the local equilibrium. The generator of $(\Theta^\epsilon(t),X^\epsilon(t))$ is $\frac{1}{\epsilon}\mathcal L+Q(x)$, acting on smooth functions of $x\in \mathbb R$ and $s_i \in S$, where 
\[\mathcal L=\left( \begin{array}{cccc}
\mathcal L_1&0&\dots&0\\
0&\mathcal L_2&\dots&0\\
&&\ddots&\\
0&0&\dots&\mathcal L_M
\end{array}\right)\]
and
\[\mathcal L_i=(s_i-x)\frac{\partial}{\partial x} + \frac{1}{2 }\frac{\partial^2}{\partial x^2}.\]
For any $x,x'\in\mathbb R$, $i,j\leq M$ and any $t\geq 0$,  we let $\Lambda_t^\epsilon$ denote the transition density function, which is the solution of the Fokker-Planck equation 

\[\left(\frac{\partial}{\partial t}- \frac{1}{\epsilon}\mathcal L^*\right)\Lambda_t^\epsilon(x,s_i|x',s_j)=\sum_\ell Q_{\ell i}(x)\Lambda_t^\epsilon(x,s_\ell|x',s_j),\]
\[\Lambda_0^\epsilon(x,s_i|x',s_j) = \delta(x-x')\delta_{ij},\]
where $\mathcal L^*$ is the adjoint of $\mathcal L$.
In particular, the transition probability over a time-step $\Delta t$ is given by $\Lambda_{\Delta t}^\epsilon$,

\[\Lambda_{\Delta t}^\epsilon(x,s_i|x',s_j)=\frac{\partial}{\partial x}\mathbb P(X^\epsilon(t+\Delta t)\leq x,\Theta^\epsilon(t+\Delta t)=s_i|X^\epsilon(t)=x',\Theta^\epsilon(t)=s_j).\]
Samples of this type of OU process with $Q(x)\equiv Q$ (i.e. $Q$ is a constant function of $x$) are shown in figure \ref{fig:OUexample}, where we that see the stationary behavior of $X^\epsilon(t)$ around $\Theta(t)$ occurs more and more rapidly as $\epsilon$ decreases.

The observation process $Y^\epsilon(t)$ depends on the hidden states since it is given by a noisy and nonlinear function of the hidden state(s). The continuous time dynamics of $Y^\epsilon(t)$ are described by the following SDE:
\[dY^\epsilon(t)=h(X^\epsilon(t))dt+dZ(t)\]
where $Z(t)$ is another independent Wiener process, and $h(\cdot)$ is some continuous and bounded function. However, \textbf{the process $\mathbf{Y^\epsilon(t)}$ is only discretely observed,} and so for some time-step $\Delta t$ there are time points $t_k=k\Delta t$ for $k=1,2,3,4,....$ such that the observations are given by a discrete sequence $Y_k^\epsilon\doteq Y^\epsilon(t_k)$, with all the observations up to time $t_k$ denoted by $Y_{1:k}^\epsilon$. The process $Y_k^\epsilon$ is given incrementally as
\[Y_{k+1}^\epsilon=Y_k^\epsilon+\int_{t_k}^{t_{k+1}}h(X^\epsilon(\tau))d\tau + \Delta Z_k\]
where $Y_{k+1}^\epsilon$ is the observation at time $t_{k+1}$ and $\Delta Z_k = Z(t_{k+1})-Z(t_k)$. The standard filtering problem for an HMM of this sort is to find an estimate of the hidden states $(\widehat\Theta^\epsilon(t_k),\widehat X^\epsilon(t_k))$ which is adapted to the filtration $\mathcal F_k^\epsilon=\sigma(Y_{1:k}^\epsilon)$, where $\sigma(Y_{1:k}^\epsilon)$ denotes the $\sigma$-algebra generated by the observations.
\subsection{The Standard Nonlinear Filter}
\label{sec:filter}
For diffusion processes, nonlinear filtering is implemented either through numerical integration or through Monte Carlo simulation (such as particle filters; see \cite{asmussen2007,cappe2005}). In most of the applied literature, it is standard to write a filtering algorithm as a recursive formula for the posterior distribution of the hidden states as new data arrives, and many filtering recursions (linear as well as nonlinear) can be derived by an application of Bayes rule. When applying Bayes rule, the essential assumptions are the Markov property and incremental independence of the noise in the observations, both of which are specified apriori in the model. These two conditions comprise what is referred to as the \textit{memoryless channel} \cite{krylov}, which allows the joint transition probabilities to be written as
\[\frac{\partial}{\partial x}\mathbb P(Y_{k+1}^\epsilon,X^\epsilon(t_{k+1})\leq x,\Theta^\epsilon(t_{k+1})=s_i|Y_{1:k}^\epsilon,X^\epsilon(t_k)=x',\Theta^\epsilon(t_k)=s_j) \]
\[= \mathbb P(Y_{k+1}^\epsilon|Y_k^\epsilon,X^\epsilon(t_{k+1})=x,X^\epsilon(t_k)=x')\times\Lambda_{\Delta t}^\epsilon(x,s_i|x',s_j).\]

Given the memoryless channel, a simple application of Bayes rule yields the recursive formulas for filtering in discrete time. The formulas for filtering in continuous time can be heuristically derived by passing the limit in $\Delta t$ in the discrete filter. In this paper the primary focus is on discretely observed processes, making it straight-forward to derive all the filtering recursions and they will always come from the Bayesian recursion.
In discrete time, the observed data is collected and retained as a sequence of real numbers from which a filtering distribution is constructed. At time $t_k$, the observed data is retained in a sequence $\{y_1,\cdots y_k\}$ that can be interpreted as the event $\{Y_{1:k}^\epsilon=y_{1:k}\}$. For any pair $s_i\in\mathcal S$ and $x\in\mathbb R$, the filtering distribution at time $t_k$ is then a function of $y_{1:k}$,
\[\pi_k^\epsilon(x,s_i) \doteq \frac{\partial}{\partial x}\mathbb P(X^\epsilon(t_k)\leq x,\Theta^\epsilon(t_k)=s_i|Y_{1:k}^\epsilon=y_{1:k}).\]

The fundamentals of Bayesian filtering are presented in the earlier works by Kushner \cite{kushner08}, Jazwinski \cite{jazwinski1970}, and Baum et al \cite{baum,baumWelch}. The Bayesian approach has since been used in many different applications. Some more recent references include Asmussen and Glynn \cite{asmussen2007}, Capp\'e, Moulines and Ryd\'en \cite{cappe2005}, Krylov et al \cite{krylov}, and Rabiner \cite{rabiner1989}. Bayesian filtering is usually done with a recursive formula, which we refer to as the \textit{Forward Baum-Welch equation}.  For any $x,x'\in\mathbb R$, $i,j\leq M$ and any $t\geq 0$, the Forward Baum-Welch equation is
\begin{equation}
\label{eq:standardFBW}
\pi_{k+1}^\epsilon(x,s_i) = \frac{1}{c_{k+1}}\sum_j\int\mathbb P(y_{k+1}|y_k,X^\epsilon(t_{k+1})=x,X^\epsilon(t_k)=x')\Lambda_{\Delta t}^\epsilon(x,s_i|x',s_j)\pi_k^\epsilon(x',s_j)dx'
\end{equation}
where $c_{k+1}$ is a normalizing constant, $\mathbb P(y_{k+1}|y_k,X^\epsilon(t_{k+1})=x,X^\epsilon(t_k)=x')$ is the likelihood function of $y_{k+1}$ given $(y_k,x,x')$, and $\pi_0(x,s_i)$ is given independent of $\epsilon$. 

The filtering recursion in (\ref{eq:standardFBW}) is difficult both to compute and to analyze because of the path-dependence in the likelihood function,
\[\mathbb P(y_{k+1}|y_k,X^\epsilon(t_{k+1})=x,X^\epsilon(t_k)=x')\]
\[ = \mathbb E\left[\mathbb P(y_{k+1}|y_k,X^\epsilon(\tau)~\forall \tau\in[t_k,t_{k+1}])|X^\epsilon(t_{k+1})=x,X^\epsilon(t_k)=x'\right]\]
\[ = C\cdot\mathbb E\left[\exp\left\{-\frac{1}{2}\left(\frac{y_{k+1}-y_k-\int_{t_k}^{t_{k+1}}h(X^\epsilon(\tau))d\tau}{\sqrt{\Delta t}}\right)^2\right\}\Bigg|X^\epsilon(t_{k+1})=x,X^\epsilon(t_k)=x'\right]\]
where $C$ is a normalizing constant (dependent on $\Delta t$). This path-dependence can be avoided by adding to the HMM an extra process that allows the filter to be written with a point-wise-dependent likelihood function. Denote this process by $V^\epsilon(t)$ and let it satisfy the differential equation
\[dV^\epsilon(t) = h(X^\epsilon(t))dt.\]
With the addition of $V^\epsilon(t)$ the observations no longer depends on a path since
\[Y_{k+1}^\epsilon=Y_k^\epsilon+V^\epsilon(t_{k+1})-V^\epsilon(t_k) + \Delta Z_k\]
We call the triplet $(V^\epsilon(t),X^\epsilon(t),\Theta^\epsilon(t))\in\mathbb R\times\mathbb R\times\mathcal S$, which is a Markov process, the augmented HMM with generator $\frac{1}{\epsilon}\mathcal L+Q(x)+h(x)\frac{\partial}{\partial v}$. For any $(v,x,s_i)\in\mathbb R\times\mathbb R\times\mathcal S$ and $(v',x',s_j)\in\mathbb R\times\mathbb R\times\mathcal S$, we let $\Gamma_t^\epsilon$ denote the Green's function of the Fokker-Planck equation given by the generator of the augmented HMM,
\[\left(\frac{\partial}{\partial t}- \frac{1}{\epsilon}\mathcal L^*-h(x)\frac{\partial}{\partial v}\right)\Gamma_t^\epsilon(v,x,s_\ell|v',x',s_j)=\sum_\ell Q_{\ell i}(x)\Gamma_t^\epsilon(v,x,s_\ell|v',x',s_j),\] 
with $\Gamma_0^\epsilon(v,x,s_\ell|v',x',s_j)=\delta(v-v')\delta(x-x')\delta_{\ell j}$.

The augmented filter is then
\[\pi_k^\epsilon(v,x,s_i) \doteq \frac{\partial^2}{\partial v\partial x}\mathbb P(V^\epsilon(t_k)\leq v,X^\epsilon(t_k)\leq x,\Theta^\epsilon(t_k)=s_i|Y_{1:k}^\epsilon=y_{1:k})\]
for which the forward Baum-Welch equation is
\[\pi_{k+1}^\epsilon(v,x,s_i) = \frac{1}{c_{k+1}}\sum_j\int\int\mathbb P(y_{k+1}|y_k,V^\epsilon(t_{k+1})=v,V^\epsilon(t_k)=v')\]
\begin{equation}
\label{eq:FBW}
\qquad\qquad\qquad\qquad\times\Gamma_{\Delta t}^\epsilon(v,x,s_i|v',x',s_j)\pi_k^\epsilon(v',x',s_j)dx'dv'
\end{equation}
where $c_{k+1}$ is another normalizing constant. The advantage of the augmented filter is that the likelihood function is now a Gaussian,
\[\mathbb P(y_{k+1}|y_k,V^\epsilon(t_{k+1})=v,V^\epsilon(t_k)=v')=C\cdot\exp\left\{-\frac{1}{2}\left(\frac{y_{k+1}-y_k-(v-v')}{\sqrt{\Delta t}}\right)^2\right\}\]
where $C$ is a normalizing constant (dependent on $\Delta t$). It can be readily verified that the marginal
posterior density $\int \pi_{k}^\epsilon(v,x,s_i) dv$ satisfies the Baum-Welch recursion (\ref{eq:standardFBW}) with the correct path-dependent likelihood function.
\subsection{Averaged Filters}
\label{sec:heuristics}
The rate at which $X^\epsilon(t)$ reverts to an approximately stationary distribution is $1/\epsilon$, and as $\epsilon\searrow 0$ we would expect there to be some stationary behavior in the filter provided that the condition in (\ref{eq:upperJumpRateBound}) holds. Indeed, this is the case and we give a direct proof in later sections, but for now we want to provide some intuition about the filter in the limit, and consider some applications. 

For any  $i,j\leq M$, let $\mu_i(x)$ be the invariant density such that $\mathcal L_i^*\mu_i(x)=0$, which can be written as
\[\mu_i(x) = \frac{1}{\sqrt\pi}e^{-(x-s_i)^2}\]
and integrates to one, $\int\mu_i(x)dx=1$, making it a probability density.
Let $\overline Q$ and $\overline h = (\overline h_{s_i})_i$ be the averages under the invariant measure,
\[\bar Q_{ij}=\int Q_{ij}(x)\mu_i(x)dx\qquad\hbox{and}\qquad\bar h_{s_i}=\int h(x)\mu_i(x)dx\qquad\hbox{for all }i\leq M.\]
Then for a fixed $T<\infty$, the process $(V^\epsilon(t),\Theta^\epsilon(t)) \in\mathbb R\times \mathcal S ,~ t\in[0,T]$
converges weakly to the process $(\overline V(t),\overline \Theta(t))\in\mathbb R\times \mathcal S,~ t\in [0,T]$
as $\epsilon\searrow 0$.\footnote{This result is proven here in Section \ref{sec:rigorousProof}; the proof follows well-known general methods (\cite{yinZhangCont} chapter 7) for weak convergence of processes.} The generator of the limit proces $(\overline V(t),\overline \Theta(t))$ is the sum of the averaged operands $\overline Q+\overline h\frac{\partial}{\partial v}$, from which it is clear that $(\overline V(t),\overline \Theta(t))$ is a Markov process. In particular, $\overline \Theta(t)$ is a continuous-time Markov chain with transition intensities $\overline Q$, and 
\[d\overline V(t) = \overline h_{\overline \Theta(t)}dt.\]
Therefore, $\overline V(t)$ is a deterministic function of the path of $\overline\Theta(t)$, and so the limiting nonlinear filter is of significantly reduced dimension and is related directly to the regime-determining path. 

Based on this result we expect the limiting nonlinear filter to be a marginal posterior probability function $\overline \pi_k$ satisfying the reduced recursion
\begin{equation}
\label{eq:averageFilter}
\overline\pi_{k+1}(s_i) = \frac{1}{\overline c_{k+1}}\sum_j\psi_{ij}(y_{k+1}|y_k) (e^{\bar Q^*\Delta t})_{ij}\overline \pi_k(s_j)
\end{equation}
where $\bar Q^*$ denotes matrix transpose of $\bar Q$, $\overline c_{k+1}$ is a normalizing constant, and $\psi$ is the likelihood function that is computed by integrating over all possible paths,
\[\psi_{ij}(y_{k+1}|y_k)=\mathbb E\left[\exp\left\{-\frac{1}{2}\left(\frac{y_{k+1}-y_k-\int_{t_k}^{t_{k+1}}\overline h_{\overline\Theta(\tau)}d\tau}{\sqrt{\Delta t}}\right)^2\right\}\Bigg|\overline \Theta(t_{k+1})=s_i,\overline \Theta(t_k)=s_j\right]. \]
The averaged filter in (\ref{eq:averageFilter}) is of significantly reduced complexity compared to the filter in (\ref{eq:standardFBW}) or the augmented filter of (\ref{eq:FBW}). There is path-dependence in the likelihood function, but this is a relatively minor complication because this dependence is on the path of a finite-state Markov chain over a relatively short interval. In contrast, the extra two dimensional state space that we need to deal with in order to filter $V^\epsilon(t_k)$ and $X^\epsilon(t_k)$ adds to the complexity of the filter. Therefore, provided that $\epsilon$ and $\Delta t$ are chosen appropriately, it is faster and almost as effective to use the average filter to compute the posterior of $\Theta^\epsilon(t_k)$. Later sections include discussion on regimes for $\epsilon$ and $\Delta t$ that make it appropriate to use this asymptotic approximation.

\subsection{Some Examples} 
\label{sec:examples}
Ornstein-Uhlenbeck-type processes with fast mean reversion are used in a variety of applications. Examples come from target tracking, biology, imaging, finance etc. 
\subsubsection{Target-Tracking Models}
Regime-switching models and nonlinear filtering are applied to target-tracking by Bar-Shalom \cite{barShalom} and also by Rozovsky \cite{petrov1999}. Suppose that a computerized-tracking system is following an object with measurements made by a video camera. For such a problem we interpret the variables as follows:
\begin{eqnarray}
\nonumber
\Theta^\epsilon(t)\qquad&=&\qquad\hbox{ a hidden maneuver to be identified,}\\
\nonumber
X^\epsilon(t)\qquad&=&\qquad\hbox{ tracked velocity}
\end{eqnarray}
where $X^\epsilon(t)$ rapidly relaxes to $\Theta^\epsilon(t)$ but in the presence of significant noise. Observations of the target are noisy measurements of the target's position,
\[Y_{k+1}^\epsilon = Y_k^\epsilon+\int_{t_k}^{t_{k+1}}h(X^\epsilon(\tau))d\tau+\Delta Z_k\]
We introduce an auxiliary field $P^\epsilon(t)$ to denote position, and a change in position of the target is $P^\epsilon(t_{k+1})-P^\epsilon(t_k) = \int_{t_k}^{t_{k+1}}h(X^\epsilon(\tau))d\tau$. Then the discrete-time observation model becomes
\[Y_{k+1}^\epsilon -Y_k^\epsilon=P^\epsilon(t_{k+1})-P^\epsilon(t_k)+\Delta Z_k\]
which can be interpreted as relative changes in measured position being equal to the relative change in the actual position plus noise. The filter for such a model is given by (\ref{eq:FBW}), and the average filter of (\ref{eq:averageFilter}) can be applied as $\epsilon\searrow 0$. Ultimately, if the computerized tracking system can positively identify $\Theta^\epsilon(t_k)$, it will be possible to have information about the sequence of maneuvers made by the target.

\textcolor{black}{The potential for an averaged filter's success in these video-tracking problems depends on the user's ability to choose $\Delta t$ such that
\[\min_{i,x}\frac{-1}{Q_{ii}(x)}\gg\Delta t\gg\epsilon.\]
Provided that $Q$ and $\epsilon$ are of different scales, the operator of the video-surveillance system should select a frame-rate for the camera that is both effective for tracking (ie. $ \min_{i,n}(-1/Q_{ii}(x))\gg\Delta t$), and consistent with the fast mean-reverting hypothesis (i.e. $\Delta t\gg\epsilon)$.
}

\subsubsection{Ornstein-Uhlebeck process for Stochastic Volatility}
\label{sec:stochasticVolExample}
Multi-scale models for stochastic volatility have been considered extensively by Fouque, Papanicolaou and Sircar \cite{FPS00} and also by Howison \cite{howison}. In this example we introduce filtering to stochastic volatility models of the form
\begin{eqnarray}
\nonumber
Y^\epsilon(t)\qquad&=&\qquad\hbox{Observed log-price}\\
\nonumber
\Theta^\epsilon(t)\qquad&=&\qquad\hbox{stochastic volatility regime} \\
\nonumber
\sqrt{h(X^\epsilon(t))}\qquad&=&\qquad\hbox{stochastic volatility}
\end{eqnarray}
Here $X^\epsilon(t)$ is the Ornstein-Uhlenbeck process defined by (\ref{eq:dX}),
the function $h(x)$ is taken to be continuous, positive, bounded and monotone, and it is assumed that there is a constant $\delta >0$ and $H< \infty$ such that
\[0<\delta\leq h(x) \leq H < \infty\qquad\forall x\in\mathbb R. \]
For example, the exponential OU model of Perell\'o, Sircar and Masoliver \cite{perello} can be modified to fit this framework. 

Letting $(X^\epsilon(t),\Theta^\epsilon(t))$ satisfy (\ref{eq:dPtheta}) and (\ref{eq:dX}), we take the logarithm of a stock price to satisfy the following SDE:
\[dY^\epsilon(t) = \left(r-\frac{1}{2}h(X^\epsilon(t))\right)dt+\sqrt{h(X^\epsilon(t))}\left(\sqrt\epsilon\rho dW(t)+\sqrt{1-\epsilon\rho^2}dZ(t)\right)\]
where $r$ is a known parameter (the interest rate), $Z(t)$ is a standard Brownian motion independent of $W(t)$, and $\rho\in(-1,0)$ models the volatility leverage effect provided that $h'(x)\geq 0 $ for all $x$,
\[\int_0^tdY^\epsilon(\tau)\cdot dh(X^\epsilon(\tau)) = \rho \int_0^th'(X^\epsilon(\tau))\sqrt{h(X^\epsilon(\tau))}d\tau.\footnote{We write ``$\int_0^tdY^\epsilon(\tau)\cdot dh(X^\epsilon(\tau))$" to denote the quadratic cross-variation $Y^\epsilon(\tau)$ and $h(X^\epsilon(\tau))$.}\]
Note that $\widetilde{W}(t) = \sqrt\epsilon\rho W(t)+\sqrt{1-\epsilon\rho^2} Z(t)$ is a standard Brownian correlated with $\frac{1}{\sqrt{\epsilon}}W$, with correlation coefficient $\rho$.
Since the log price $Y^\epsilon(t)$ is observed discretely (e.g. daily), we write the observation model as follows:
\[Y_{k+1}^\epsilon-Y_k^\epsilon =r\Delta t-\frac{1}{2}\int_{t_k}^{t_{k+1}}h(X^\epsilon(\tau))d\tau+\sqrt\epsilon\rho\int_{t_k}^{t_{k+1}}\sqrt{h(X^\epsilon(\tau))}dW(\tau)+\sqrt{1-\epsilon\rho^2}\int_{t_k}^{t_{k+1}}\sqrt{h(X^\epsilon(\tau))}dZ(\tau)\]
\[=_d r \Delta t-\frac{1}{2}\int_{t_k}^{t_{k+1}}h(X^\epsilon(\tau))d\tau+\sqrt\epsilon\rho\int_{t_k}^{t_{k+1}}\sqrt{h(X^\epsilon(\tau))}dW(\tau)
+\sqrt{(1-\epsilon\rho^2)\int_{t_k}^{t_{k+1}}h(X^\epsilon(\tau))d\tau}\cdot\mathcal Z_k\]
where the equality is in distribution and $\mathcal Z_k\sim iid N(0,1)$. 

As before, we introduce two auxiliary fields, namely $(V^\epsilon(t),\zeta^\epsilon(t))$, where
\[dV^\epsilon(t) = (1-\epsilon\rho^2)h(X^\epsilon(t))dt\qquad\hbox{and}\qquad d\zeta^\epsilon(t)=\sqrt\epsilon\rho\sqrt{h(X^\epsilon(t))}dW(t)\]
with initial conditions $V^\epsilon(0) = 0=\zeta^\epsilon(0)$, for which we see that $(V^\epsilon(\cdot),\Theta^\epsilon(\cdot))\Rightarrow (\overline V(\cdot),\overline\Theta(\cdot))$ as $\epsilon\searrow 0$, and $\zeta^\epsilon(t)\rightarrow 0 $ a.s. as $\epsilon\searrow 0$ for all $t<\infty$. Applying these limits, we have the following weak limit of the observations,
\[Y_{k+1}^\epsilon-Y_k^\epsilon\Rightarrow \mu\Delta t-\frac{1}{2}\int_{t_k}^{t_{k+1}}\bar h_{\bar\Theta(\tau)}d\tau+\sqrt{\int_{t_k}^{t_{k+1}}\bar h_{\bar\Theta(\tau)}d\tau}\cdot\mathcal Z_k\]
and by the same token that shows $V^\epsilon(\cdot)$ converges, we have the following averaging of the leverage effect,
\[\int_{t_k}^{t_{k+1}}dY^\epsilon(\tau)\cdot dh(X^\epsilon(\tau))\Rightarrow\rho \int_{t_k}^{t_{k+1}} (\overline{h'\sqrt h})_{\overline \Theta(\tau)}d\tau\]
where $(\overline{h'\sqrt h})_{s_i} = \int h'(x)\sqrt{h(x)}\mu_i(x)dx$, both limits taken as $\epsilon\searrow 0$. 

The standard nonlinear filter for this problem is more complicated than the one given in equation (\ref{eq:FBW}), but we will see later on that a simple generalization of theorem \ref{thm:averagedFilter} of Section \ref{sec:rigorousProof} leads to an averaged filter as $\epsilon\searrow 0$,
\[\overline \pi_{k+1}(s_i) = \frac{1}{\overline c_{k+1}}\sum_j\psi_{ij}(y_{k+1}|y_k)\left(e^{\bar Q^*\Delta t}\right)_{ij}\overline \pi_k(s_j)\]
where the likelihood function has the form
\[\psi_{ij}(y_{k+1}|y_k)=\mathbb E\left[\frac{\exp\left\{-\frac{1}{2}\left(\frac{y_{k+1}-y_k-\mu\Delta t+\frac{1}{2}\int_{t_k}^{t_{k+1}}\overline h_{\overline \Theta(\tau)}d\tau}{\sqrt{\int_{t_k}^{t_{k+1}}\overline h_{\overline \Theta(\tau)}d\tau}}\right)^2 \right\}}{\sqrt{\int_{t_k}^{t_{k+1}}\overline h_{\overline \Theta(\tau)}d\tau}}\Bigg| \overline \Theta(t_{k+1})=s_i,\overline \Theta(t_k)=s_j\right]\]
with $\overline h_{s_i} = \int h(x)\mu_i(x)dx=\frac{1}{\sqrt\pi}\int h(x)e^{-(x-s_i)^2}dx$. We will revisit this example in Section \ref{sec:rigorousProof} after we have gone over the rigorous proof of the averaged filter, at which point we will be able to give a detailed explanation as to why this stochastic volatility filter has the asymptotic averaging that we show here. From the point of view stochastic volatility modeling, we see that nonlinear filters for tracking regime changes in the volatility can be significantly simplified if we have fast mean reversion and if the observation interval is chosen appropriately.

\textcolor{black}{The potential for the averaged filter to be effective in tracking the state of volatility will depend on the data. Time scales in financial data have been observed, with fast mean-reversion happening on the order of 1-2 days (or less), and slow mean-reversion happening on the order of weeks and months. Therefore, a choice of $\Delta t$ corresponding to a few days or a week would allow the averaged filter to track regime-changes effectively. It should also be mentioned that there are time-dependent effects in financial data (e.g. daily intra-day trading patterns) that should be included in a stochastic volatility model, but the averaged filter has the ability to track a time-independent regime-process without estimating/tracking the time-dependent components (see \cite{FPSS04}). }


\subsection{Issues with $\mathbf\epsilon$ and $\mathbf{\Delta t}$}
\label{sec:issues}
There are some subtle issues in fast mean reversion models and filtering that can be easily overlooked. In this subsection we highlight some of them that can arise when seemingly small changes are made to the proposed model. They are from a general class of fast mean reverting models that may look similar to the model in this paper, but are in fact quite different. In particular their limiting filters are not the same.

\subsubsection{Allowing for Relaxation Between Observations}
For a fixed $\epsilon>0$, the choice to implement the filter in (\ref{eq:averageFilter}) is appropriate provided that the sampling rate of observations, $\Delta t$, is chosen correctly. If $\Delta t$ is too small (i.e. $\Delta t<<\epsilon$) then there will not be enough time between observations for the diffusion distribution of $X^\epsilon$ to relax. If the average filter is applied in such a situation, there will be estimation error that is not present when the optimal nonlinear filter is used. Therefore, it is important to impose a condition such as $\Delta t>\epsilon$ to decide if it is appropriate to use the average filter in (\ref{eq:averageFilter}). In other words, such a condition ensures that the average time to mean reversion, $\epsilon$, in the diffusion process is faster than the rate at which new observations arrive. The interval between observations, $\Delta t$, must also not be too large because then the filter will be neglectful of regime changes. In applications where there is regime change and fast mean reversion, and the time scales for these two effects are well separated, the observation time should be chosen to be intermediate between the two,
\begin{equation}
\label{eq:timescales}
\epsilon \ll \Delta t \ll \min_{i,x}\frac{-1}{Q_{ii}(x)}.
\end{equation}

The above discussion suggests that fast mean-reversion asymptotics cannot be applied to the Zakai equation \cite{rozovsky1991,zakai1969}, which gives the continuous time filter.  This is because averaged filters only work when there is some coarse-graining of the observations (for example, by discretizing the observation time) and when the rate of mean reversion is faster than the observation sampling rate.

\subsubsection{Different  Scalings}
\label{sec:papanicolaou2010}
In this section we present an example to illustrate some subtleties in the scaling of stochastic differential equations. The example is the focus of the paper by Papanicolaou \cite{papanicolaou2010} in which there is a limiting filter of reduced dimension, but the type of averaging that is used to obtain it is different and thus the limiting filter is considerably different from that in (\ref{eq:averageFilter}). 

We consider the following video-tracking model:
\begin{eqnarray}
\nonumber
\Theta(t)\qquad&=&\qquad\hbox{target velocity}\\
\nonumber
X^\epsilon(t)\qquad&=&\qquad\hbox{camera velocity}\\
\nonumber
P^\epsilon(t)\qquad&=&\qquad\hbox{position error}\\
\nonumber
Y_k^\epsilon\qquad&=&\qquad\hbox{noisy measurement of position error}
\end{eqnarray}
Suppose the dynamics of $\Theta(t)$ have intensity matrix $Q$ which no longer depends on $X^\epsilon(t)$, and the diffusion is no longer scaled by $1/\sqrt\epsilon$,
\begin{eqnarray}
\nonumber\frac{\partial}{\partial t}\mathbb P(\Theta(t)=s_i) &=& \sum_jQ_{ji}\mathbb P(\Theta(t)=s_j)\\
\nonumber
dX^\epsilon(t) &=& \frac{1}{\epsilon}(\Theta(t)-X^\epsilon(t))dt+ dW(t)\\
\nonumber
dP^\epsilon(t) &=& \frac{1}{\epsilon}(\Theta(t)-X^\epsilon(t))dt\\
\nonumber
Y_k^\epsilon &=& h\left(P^\epsilon(t_k)\right)+Z_k
\end{eqnarray}
where $Z_k\sim N(0,1)$. The main differences between the model here and in the rest of the paper is the lack of a scaling factor before the $dW_t$ term and the addition of a scaling term in the dynamics of the auxiliary field (in this case we have $P^\epsilon$ as the auxiliary field). In \cite{papanicolaou2010}, the limit theory used to obtain the average filter for this model is simpler than the theory presented in the rest of this paper, but the asymptotic average of the filter is a of an entirely different nature.

The filter for this model is
\[\pi_k^\epsilon(p,x,s_i) = \frac{\partial}{\partial x}\frac{\partial }{\partial p}\mathbb P(P^\epsilon(t_k)\leq p,X^\epsilon(t_k)\leq x,\Theta(t_k)=s_i|Y_{1:k}^\epsilon=y_{1:k})\]
and satisfies the following forward Baum-Welch equation:
\[\pi_{k+1}^\epsilon(p,x,s_i)=\frac{1}{c_{k+1}}e^{-\frac{1}{2}(y_{k+1}-h(p))^2}\sum_j\int\int \Lambda_{\Delta t}^\epsilon(p,x,s_i|p',x',s_j)\pi_k^\epsilon(p',x',s_j)dx'dp'\]
where $\Lambda_{\Delta t}^\epsilon(~\cdot~)$ is the joint transition kernel of the hidden process. By a generalization of the Kramers-Smoluchowski theorem \cite{papanicolaou2010,freidlin04}, we have
\[X^\epsilon(t)\rightarrow \Theta(t)\qquad\qquad\hbox{a.s. pointwise in time},\]
from which it follows that
\[P^\epsilon(t_{k+1})-P^\epsilon(t_k)\rightarrow\Theta(t_{k+1})-\Theta(t_k)-\Delta W_k\qquad\hbox{a.s. pointwise in time}.\]
These limits allow us to obtain the following averaged filter in the limit
\[\overline \pi_{k+1}(p,s_i) = \frac{1}{\overline c_{k+1}}e^{-\frac{1}{2}(y_{k+1}-h(p))^2}\sum_j\int e^{-\frac{1}{2}\left(\frac{p-p'-(s_i-s_j)}{\sqrt{\Delta t}}\right)^2}\left(e^{Q^*\Delta t}\right)_{ij}\overline\pi_k(p',s_j)dp'\]

This filter is of reduced dimension, but it is not entirely reduced to a filter for the state $\Theta(t_k)$. In fact, such a fully reduced filter is not possible because the likelihood function has not been averaged. Therefore, although the theory used to obtain the limit is simpler and the convergence stronger, the limiting filter is nevertheless more complicated than the averaged filter of (\ref{eq:averageFilter}).

\section{Rigorous Averaging For Filter}
\label{sec:rigorousProof}
In this section we reintroduce the model in a probabilistic setting, which will allow us to rigorously prove the limit theorem, giving us the average filter of (\ref{eq:averageFilter}). The statements and proofs that build the theory for the averaged filter use pathwise weak convergence and test functions, and so it is necessary to redefine all the dynamics of the processes in their weak or integrated form. In particular, convergence of the solutions to a martingale problem is the essential tool in showing that the processes $(\Theta^\epsilon(t),V^\epsilon(t))~,0\leq t\leq T,$ converge path-wise in the weak sense.

\subsection{Preparation for Theorem}
\label{sec:reformulation}
Consider the Markov process $(\Theta^\epsilon(t),X^\epsilon(t),V^\epsilon(t))$ for $\epsilon>0$ and $t\in[0,T]$, on the probability space $(\Omega,(\mathcal F_t)_{t\leq T},\mathbb P)$. Let $\Theta^\epsilon(t)$ be a jump process taking value in the space of finite values $\mathcal S=\{s_1,.....,s_M\}$. The process $X^\epsilon(t)$ is an Ornstein-Uhlenbeck (OU) process with a shifting-mean,
\begin{equation}
\label{eq:OUprocess}
dX^\epsilon(t) = \frac{1}{\epsilon}\left(\Theta^\epsilon(t)-X^\epsilon(t)\right)dt+\frac{1}{\sqrt{\epsilon}}dW(t)
\end{equation}
where $W(t)$ is an independent, standard Brownian motion and $\epsilon>0$ is an arbitrarily small parameter. From this, the third process is defined as 
\[V^\epsilon(t)= V_0+\int_0^th(X^\epsilon(\tau))d\tau\]
where $h:\mathbb R\rightarrow\mathbb R$ is bounded and continuous.

The space $\Omega$ can be defined as
\[\Omega = D\left([0,T];\mathcal S\right)\times C\left([0,T];\mathbb R^2\right)\]
were $D\left([0,T];\mathcal S\right)$ can be taken to be here the set of right-hand continuous piece-wise constant functions in $\mathcal S$, 
$C\left([0,T];\mathbb R^2\right)$ are the continuous functions in $\mathbb R^2$, and $\mathcal F_t,~ 0\leq t\leq T$ is the associated filtration. There is a matrix operator $Q$  where 
\[Q(x)=(Q_{ji}(x))\in\mathbb R^{M\times M},\]
which is a transition rate matrix for all $x\in\mathbb R$ such that for any bounded function $g(\theta):\mathcal S\rightarrow \mathbb R$, 
\[Q(x)g(\theta)\Big|_{\theta=s_i} = \sum_jQ_{ij}(x)g(s_j),\]
 We assume the conditions set forth in (\ref{eq:upperJumpRateBound}) and (\ref{eq:lowerJumpRateBound}) of Section \ref{sec:formulation} hold (with constants $0<\alpha\leq\beta<\infty$) so that $\Theta^\epsilon(t)$ communicates on all possible states of $\mathcal S$, that jumps cannot occur at an infinite rate, and that there are no cemetery states. We also take the generators $\mathcal L_i$ to be the same as they were in section (\ref{sec:formulation}), and we take the matrix of these operators $\mathcal L$ to be the same as well. Then for any scalar-valued function $g(s,x,v)$ with compact support, two derivatives in $x$ and one in $v$ (both bounded), the dynamics of the Markov process $(\Theta^\epsilon(t),X^\epsilon(t),V^\epsilon(t))$ are determined by the martingale,
\[\mathbb E\left[g(\Theta^\epsilon(t),X^\epsilon(t),V^\epsilon(t))-\int_0^t\left(\frac{1}{\epsilon}\mathcal L+Q(X^\epsilon(\tau))+h(X^\epsilon(\tau))\frac{\partial}{\partial v}\right)g(\Theta^\epsilon(\tau),X^\epsilon(\tau),V^\epsilon(\tau))d\tau\Big|\mathcal F_s\right]\]
\begin{equation}
\label{eq:martingaleForm}
=g(\Theta^\epsilon(s),X^\epsilon(s),V^\epsilon(s))
-\int_0^s\left(\frac{1}{\epsilon}\mathcal L+Q(X^\epsilon(\tau))+h(X^\epsilon(\tau))\frac{\partial}{\partial v}\right)g(\Theta^\epsilon(\tau),X^\epsilon(\tau),V^\epsilon(\tau))d\tau
\end{equation}
holding for any test function $g$ and $0\leq s < t\leq T$. Equation (\ref{eq:martingaleForm}) is a weak form of the forward equation for the Markov process $(\Theta^\epsilon(t),X^\epsilon(t),V^\epsilon(t))$.

Observations on $(\Theta^\epsilon(t),X^\epsilon(t),V^\epsilon(t))$ are given at the discrete times $t_k=k\Delta t$ by the process $Y_k^\epsilon$,
\[Y_{k+1}^\epsilon=Y_k^\epsilon+V^\epsilon(t_{k+1})-V^\epsilon(t_k)+\Delta Z_k\]
where $\Delta Z_k\sim iid N(0,\Delta t)$. 

\subsection{Theorem for Averaged Filter}
\label{sec:theorem}
We consider a sequence $y_{1:k}=\{y_1,\dots,y_k\}$ to be the observed realization of $Y_{1:k}^\epsilon$, for which we have the following asymptotic theory for the filter.

\begin{theorem}
\label{thm:averagedFilter}
Assume (\ref{eq:upperJumpRateBound}) is true. Then for any bounded function $g:\mathcal S\rightarrow \mathbb R$, point-wise for any $y_{1:k}=\{y_1,y_2,...,y_k\}\in\mathbb R^k$ the optimal filter has the following limit as $\epsilon\searrow 0$:  
\[\mathbb E\left[g(\Theta^\epsilon(t_k))\Big|Y_{1:k}^\epsilon=y_{1:k}\right]\rightarrow \sum_ig(s_i)\bar\pi_k(s_i)\qquad\hbox{ as }\epsilon\searrow 0\]
where 
\begin{equation}
\label{eq:ergodicFBW}
\bar\pi_{k+1}(s_i) = \frac{1}{\bar c_{k+1}}\sum_j\psi_{ij}(y_{k+1}|y_k)\left(e^{\bar Q^*\Delta t}\right)_{ij}\bar\pi_k(s_j)
\end{equation}
with 
\[
\psi_{ij}(y_{k+1}|y_k) = \mathbb E\left[\exp\left\{-\frac{1}{2\Delta t}\left(y_{k+1}-y_k-\int_{t_k}^{t_{k+1}}\bar h_{\overline\Theta(\tau)}d\tau\right)^2\right\}\Big|\overline\Theta(t_{k+1}) = s_i,\overline \Theta(t_k)=s_j\right]
\] 
and $\bar c_{k+1}$ being a normalizing constant. 
\end{theorem}

One main component of the proof of theorem \ref{thm:averagedFilter} is the path-wise weak limit of $(\Theta^\epsilon(\cdot),V^\epsilon(\cdot))$,
\begin{equation}\label{eq:weakLimit}
\left(\Theta^\epsilon(\cdot),V^\epsilon(\cdot)\right)\Rightarrow (\overline\Theta(\cdot),\overline V(\cdot))\end{equation}
where $(\overline\Theta(\cdot),\overline V(\cdot))$ is a Markov process with generator $\bar Q+diag(\bar h )\frac{\partial}{\partial v}$. Appendix \ref{sec:details} goes through the steps proving that the limit in (\ref{eq:weakLimit}) holds given the setup of Section \ref{sec:reformulation}. 

The other component of the proof of theorem \ref{thm:averagedFilter} is the use of the unnormalized posterior. The kernel used to obtain the unnormalized posterior is a likelihood function, and so in spirit it is similar to the change of measure used in the Zakai equation (in the case of continuous-time observations) as was done by Zakai \cite{zakai1969}, and also is also a part of Rozovsky's proof of the uniqueness of solutions to the Zakai equation \cite{rozovsky1991}. 

\noindent\textbf{Proof of Theorem \ref{thm:averagedFilter}:} Without loss of generality, let $\Delta t = 1$. 
For any vector $(y_1,...,y_k)\in\mathbb R^k$ define the likelihood function of $V^\epsilon(\cdot)$,
\[\mathcal M_k(y_{1:k};V^\epsilon(\cdot)) = \prod_{i=1}^{k-1}\exp\left\{-\frac{1}{2}\left(y_{i+1}-y_i-(V^\epsilon(t_{i+1})-V^\epsilon(t_i))\right)^2\right\}.\]
The posterior expectation can be written as a function of the observed data in terms of $\mathcal M_k$,
\[\mathbb E\left[g(\Theta^\epsilon(t_k))\Big|Y_{1:k}^\epsilon=y_{1:k}\right]=\frac{\mathbb E\left[\mathcal M_k(y_{1:k};V^\epsilon(\cdot))g(\Theta^\epsilon(t_k))\right] }{\mathbb E\left[ \mathcal M_k(y_{1:k};V^\epsilon(\cdot))\right]} \]
We can now apply the weak convergence theorem to the ratio on the right. From the weak limit of $(\Theta^\epsilon(\cdot),V^\epsilon(\cdot))$ in (\ref{eq:weakLimit}) this is equal to a function of $y_{1:k}$ which converges point-wise,
\[\frac{\mathbb E\left[\mathcal M_k(y_{1:k};V^\epsilon(\cdot))  g(\Theta^\epsilon(t_k))\right] }{\mathbb E\left[ \mathcal M_k(y_{1:k};V^\epsilon(\cdot)) \right]}\rightarrow\frac{\mathbb E\left[\mathcal M_k(y_{1:k};\overline V(\cdot))g(\overline \Theta(t_k))\right] }{\mathbb E\left[ \mathcal M_k(y_{1:k};\overline V(\cdot))\right]} \]
as $\epsilon\searrow 0$. Letting $ \overline Y_k$ denote the process
\[\overline Y_k =\int_{t_k}^{t_{k+1}}\bar h_{\overline \Theta(\tau)}d\tau +Z_k\]
it is easily seen that $\mathcal M_k(y_{1:k};\overline V(\cdot))$ is the likelihood function associated with the filtering problem as $\epsilon\searrow 0$. Therefore, the limit ratio can be identified with a conditional expectation and we have
\[\mathbb E\left[g(\Theta^\epsilon(t_k))\Big|Y_{1:k}^\epsilon=y_{1:k}\right]\rightarrow \mathbb E\left[g(\overline\Theta(t_k))\Big|\overline Y_{1:k}=y_{1:k}\right]\]
as $\epsilon\searrow 0$ pointwise for $y_{1:k}$.

In connection with the Baum-Welch equation, the posterior expectation is written as
\[\mathbb E\left[g(\overline\Theta(t_k))\Big|\overline Y_{1:k}=y_{1:k}\right]=\sum_ig(s_i)\bar\pi_k(s_i)\]
where the posterior distribution $\bar \pi_k$ is given recursively as 
\[\bar\pi_k(s_i) = \frac{1}{c_k}\sum_j\frac{\partial}{\partial y}\mathbb P(\overline Y_k\leq y|\overline Y_{k-1}=y_{k-1}\overline\Theta(t_k)=s_i,\overline\Theta(t_{k-1})=s_j)\Bigg|_{y=y_k}\left(e^{\bar Q^*\Delta t_k}\right)_{ji}\bar\pi_{k-1}(s_j)\]

\[= \frac{1}{c_k}\sum_j\psi_{ji}(y_k|y_{k-1})\left(e^{\bar Q^*\Delta t_k}\right)_{ji}\bar\pi_{k-1}(s_j).\]
This completes the proof of the theorem. $\blacksquare$\\
\subsection{Stochastic Volatility Example (Revisited)}
We now provide some details to show why theorem \ref{thm:averagedFilter} applies to the stochastic volatility model of Section \ref{sec:stochasticVolExample}. Recall that the example had two hidden states
\begin{eqnarray}
\nonumber
\Theta^\epsilon(t)&=&\hbox{the hidden volatility regime}\\
\nonumber
X^\epsilon(t) &= &\hbox{the hidden volatility process}
\end{eqnarray}
with two auxiliary fields
\[dV^\epsilon(t) = (1-\epsilon\rho^2)h(X^\epsilon(t))dt\qquad\hbox{and}\qquad d\zeta^\epsilon(t)=\sqrt\epsilon\rho\sqrt{h(X^\epsilon(t))}dW(t)\]
where $V^\epsilon(0)=0=\zeta^\epsilon(0)$, and observations were given by
\[Y_{k+1}^\epsilon-Y_k^\epsilon =_d \left(\mu\Delta t-\frac{1}{2}\left(V^\epsilon(t_{k+1})-V^\epsilon(t_k)\right)\right)+\zeta^\epsilon(t_{k+1})-\zeta^\epsilon(t_k) - \sqrt{V^\epsilon(t_{k+1})-V^\epsilon(t_k)}\cdot\mathcal Z_k\]
where $\mathcal Z_k\sim iid N(0,1)$. It was also mentioned in Section \ref{sec:stochasticVolExample} that path-wise convergence holds for $\Theta^\epsilon(\cdot)$ and the auxiliary fields,
\[(\zeta^\epsilon(\cdot),V^\epsilon(\cdot),\Theta^\epsilon(\cdot))\Rightarrow (0,\overline V(\cdot),\overline \Theta(\cdot))\]
as $\epsilon\searrow 0$.

Now, in order to apply the steps in the proof of theorem \ref{thm:averagedFilter}, the likelihood function is defined as
\[\mathcal M_k(y_{1:k};V^\epsilon(\cdot),\zeta^\epsilon(\cdot)) = \prod_{i=1}^{k-1}\frac{\exp\left\{-\frac{1}{2}\left(\frac{y_{i+1}-y_i-\mu\Delta t+\frac{1}{2}(V^\epsilon(t_{i+1})-V^\epsilon(t_i))-\zeta^\epsilon(t_{i+1})-\zeta^\epsilon(t_i)}{\sqrt{(V^\epsilon(t_{i+1})-V^\epsilon(t_i))\Delta t}}\right)^2\right\}}{(V^\epsilon(t_{i+1})-V^\epsilon(t_i))\Delta t},  \]
and because $0<\delta \leq h(x) \leq H < \infty$, we have $ \mathcal M_k(y_{1:k};V^\epsilon(\cdot),\zeta^\epsilon(\cdot)) $ is bounded for all $t_k<\infty$. 
Then, using the weak convergence we have that for any bounded function $g:\mathcal S\rightarrow\mathbb R$ 
\[\mathbb E[g(\Theta^\epsilon(t_k))|Y_{1:k}^\epsilon=y_{1:k}]\rightarrow \sum_ig(s_i)\bar\pi_k(s_i)\]
as $\epsilon\searrow 0$, where $\bar \pi_k$ is given by theorem \ref{thm:averagedFilter} but with 
\[\psi_{ij}(y_{i+1}|y_i) = \mathbb E\left[\frac{\exp\left\{-\frac{1}{2}\left(\frac{y_{i+1}-y_i-\mu\Delta t+\frac{1}{2}\int_{t_i}^{t_{i+1}}\bar h_{\bar\Theta(\tau)}d\tau}
{\sqrt{\left(\int_{t_i}^{t_{i+1}}\bar h_{\bar\Theta(\tau)}d\tau\right)\Delta t}}\right)^2\right\}}
{\sqrt{\left(\int_{t_i}^{t_{i+1}}\bar h_{\bar\Theta(\tau)}d\tau\right)\Delta t} } \right]\]
for all $i\leq k-1$.

\textcolor{black}{One disadvantage to this averaged filter is that $\rho$ does not appear. Leverage is present for all $\epsilon>0$, but depends on the path of $X^\epsilon$. Therefore, the averaged filter is useful only in settings where volatility's regime appears independent of returns. For instance, suppose option prices have independent short-term dynamics, but the common level of volatility remains constant for extended periods and does not have any short-term effects on returns. If leverage is a considerable concern, it can make in an impact on the averaged filter if the model includes a slowly-scaled OU which depends on $\Theta$ (see Fouque, Papanicolaou and Sircar \cite{FPS00}). Another option is to add jumps in a manner that is similar Bakshi, Cao and Chen \cite{BCC97},
\[dY^\epsilon(t)=\left(\mu-\frac{1}{2}h(X^\epsilon(t))\right)dt+\sqrt{h(X^\epsilon(t))}dW(t)-J(t)d\Theta^\epsilon(t)-\nu dt\]
where $d\Theta^\epsilon(t)=\Theta^\epsilon(t)-\Theta^\epsilon(t^-)$, $J(t)$ is a random variable that depends on the sign of $d\Theta^\epsilon(t)$, and $\nu$ is a compensator.}
\section{Numerical Simulations}
\label{sec:numerics}
\textcolor{black}{In this section we carry out some numerical simulations of the proposed model and implement the filtering algorithms. Methods for simulating the HMM and computing the filter are described in detail, and simulated results verify some basic properties of the averaged filter. In particular, Sections \ref{sec:simulation} and \ref{sec:computeFilter} describes the numerical methods, and Sections \ref{sec:filterError} and \ref{sec:empirical} explore how the filter performs for various $\epsilon$ and $\Delta t$.}

\textcolor{black}{In estimating the hidden states, it is optimal to use the filter from Section \ref{sec:filter}, referred to as the `optimal filter', but Theorem \ref{thm:averagedFilter} proved that the averaged filter is equivalent as $\epsilon\searrow 0$.  Section \ref{sec:filterError} will focus on how the averaged filter compares to the optimal filter, and for this reason the simulations use a simplified linear model under which it is easier to compute the optimal filter. The main conclusions are 
\begin{enumerate}
\item that the averaged filter is indeed optimal as $\epsilon\searrow 0$, 
\item that $\Delta t$ needs to be chosen within a certain range depending on the mean time of reversion and the expected switching times of the hidden Markov chain, as in (\ref{eq:timescales}).
\end{enumerate}}

\subsection{Generating the Simulation Data}
\label{sec:simulation}
The simulation occurs over a time interval $[0,T]$, and given time-step $\Delta t$ we will have observations $Y_k^\epsilon$ occurring at a set of equidistant time-points $\{t_k\}_{k=0}^N$

\[0=t_0<t_1<t_2<....<t_N=T\]
with $t_{k+1}=t_k+\Delta t$. Given the linear observation function $h(x)=h\cdot x$, the observations are given simply by
\[Y_{k+1}^\epsilon=Y_k^\epsilon+h\cdot\int_{t_k}^{t_{k+1}}X^\epsilon(s)ds+\Delta Z_k\]
where $\Delta Z_k =Z(t_{k+1})-Z(t_k)$. The task of simulating consists of two parts: (1) how to generate the hidden states, and (2) how to generate the observations.

We'll assume the dynamics of $\Theta(t)$ are homogeneous (i.e. they are not affected by $X^\epsilon(t)$; $Q$ does not have an `$x$' argument). To simulate the processes, we need to selected a finer time-step, which we will call $\widetilde\Delta t$. For some positive integer $m$, we define $\widetilde\Delta t$ as
\[\widetilde\Delta t=\frac{1}{m}\Delta t,\qquad\hbox{for some }m\in\mathbb Z^+.\]
Thus, the observations are $m$-many $\widetilde\Delta t$'s apart, and the simulated state processes are generated at the $m$-many points in between,
\[\widetilde t_{\ell+1}=\widetilde t_\ell+\widetilde \Delta t\quad\hbox{for }\ell=0,1,2,\dots,m\cdot N-1.\]
At each time $\widetilde t_\ell$, a discrete-time Markov process $(\widetilde\Theta_\ell,\widetilde X_\ell^\epsilon)$ is sampled sequentially according to Algorithm \ref{alg:samples}.

\begin{algorithm}
\caption{Sequential Sampling of State-Processes}
\label{alg:samples}
\begin{algorithmic}
 \STATE \textbf{1.} set $a=\exp(-\widetilde\Delta t/\epsilon)$.
\STATE \textbf{2.} set the matrix $p=[p_{ij}]$ where $p_{ij}\doteq\left\{\exp\{\widetilde\Delta tQ^*\}\right\}_{ij}$
\STATE \textbf{3.} Sequentially sample the states:
\FOR{$\ell = 0 \to m\cdot N-1$}
\STATE\textbf{a.} $ \widetilde \Theta_{\ell+1}\gets \theta$ where $\theta$ is sampled from $p$ conditional on $\widetilde\Theta_\ell$,
\STATE \textbf{b.} $ \widetilde X^\epsilon_{\ell+1}\gets a\widetilde X^\epsilon_\ell+(1-a)\widetilde\Theta_\ell + \sqrt{\frac{1-a^2}{2}}\mathcal W_\ell$ where $\mathcal W_\ell\sim iid N(0,1).$
\ENDFOR
\end{algorithmic}
\end{algorithm}
Given the discrete process from Algorithm \ref{alg:samples}, we approximate $(\Theta(\cdot),X^\epsilon(\cdot))$ with a piece-wise constant process $(\widetilde\Theta(\cdot),\widetilde X^\epsilon(\cdot))$, defined as follows:\\

\begin{equation}
\label{eq:pointProcess}
(\widetilde\Theta(t),\widetilde X^\epsilon(t)) \doteq \sum_{\ell = 0}^{m\cdot N-1} (\widetilde\Theta_\ell,\widetilde X^\epsilon_\ell )\cdot\mathbf 1_{\{t\in [\widetilde t_\ell,\widetilde t_{\ell+1})\}}\qquad\hbox{for }t\in[0,T].
\end{equation}
The process in (\ref{eq:pointProcess}) converges path-wise in the weak sense to $(\Theta(\cdot),X^\epsilon(\cdot))$ as $m\rightarrow \infty$ (the proof of this convergence is similar to the convergence proofs in Appendix \ref{app:convergenceMartingale}, and is also covered in \cite{ethierKurtz}).

Finally, we approximate the observations with Riemann sums,
\begin{equation}
\label{eq:observedApprox}
\widetilde Y_{k+1}^\epsilon=\widetilde Y_k^\epsilon+\widetilde \Delta t ~h\cdot \sum_{\ell=mk+1}^{m(k+1)}\widetilde X_\ell^\epsilon+\Delta Z_k\qquad\hbox{for }k=0,1,2,\dots,N-1.
\end{equation}
Indeed, it can be further proved that $\widetilde Y_k^\epsilon \Rightarrow Y_k^\epsilon$ as $m\rightarrow\infty$, and so our simulations will assume that $\widetilde Y^\epsilon=Y^\epsilon$ and choose $m$ so that the Riemann sum is close to an integral.
 
 \subsection{Computing the Filter}
 \label{sec:computeFilter}
Suppose that we have used Algorithm \ref{alg:samples} and (\ref{eq:observedApprox}) to generate the observed the process $Y_{1:k}^\epsilon=y_{1:k}$. Algorithm \ref{alg:samples} is used again to generate samples for a particle filter (see \cite{asmussen2007,cappe2005}) or a Rao-Blackwellized filter (see \cite{gustafsson2002,gustafsson2005}), 
\[\{\widetilde\theta_\ell^{(n)},\widetilde X_\ell^{\epsilon,{(n)}}\}_{\ell=0}^{m\cdot N}\qquad\hbox{for }n=1,2,3,\dots,R\]
where $R$ is the number of samples generated. The particle filter consists of a set of weights $\{\omega_k^{(n)}\}_{n,k}$ such that for $R$ large we have the following approximation:
\[\mathbb E\left[g(X^\epsilon(t_k),\Theta(t_k))\Big|Y_{1:k}^\epsilon=y_{1:k}\right]\approx \sum_{n=1}^Rg(\widetilde\Theta_{km}^{(n)},\widetilde X_{km}^{\epsilon,(n)})\cdot\omega_k^{(n)}\qquad\hbox{for }k=0,1,2,\dots,N.\]
where the weights are computed as
\[\omega_{k+1}^{(n)}=\frac{1}{c_{k+1}}\exp\left\{-\frac{1}{2}\left(\frac{y_{k+1}-y_k-\widetilde \Delta t~h\cdot\sum_{\ell=km+1}^{m(k+1)} \widetilde X_\ell^{\epsilon,(n)}}{\sqrt{\Delta t}}\right)^2\right\}\cdot\omega_k^{(n)}\]
with $c_{k+1}$ being a normalizing constant. This algorithm of sampling and weighting is usually accompanied by sampling-importance-resampling (SIR), which basically replaces the samples with boot-strap samples if relative importance starts to favor a select-few particles. For instance if the number of important particles falls below some critical level we should perform SIR,
\[1/\sum_n(\omega_k^{(n)})^2\leq \eta R,\qquad\hbox{then perform SIR}\]
where $\eta\in[0,1]$ is the required fraction of important particles. SIR involves the following: (1) for each $n$ we sample a random variable $V_n$ from $\left\{(\widetilde\Theta_{m(k+1)}^{(n)},\widetilde X_{m(k+1)}^{\epsilon,(n)})\right\}_{n=1}^R$ with probabilities $\left\{\omega_k^{(n)}\right\}_{n=1}^R$; (2) after $R$-many $V$'s have been obtained, set
\[(\widetilde\Theta_{mk+\ell+1}^{(n)},\widetilde X_{mk+\ell+1}^{\epsilon,(n)})=V_n,\qquad\omega_k^{(n)}=1/R\qquad\hbox{for each }n.\]
Basically, SIR is a variance reduction technique that will prevent the particle filter from relying on relatively few samples. Algorithm \ref{alg:particleFilter} summarizes the steps taken to recursively to compute the particle filter upon the arrival of the $(k+1)^{th}$ observation. 

\begin{algorithm}
\caption{Recursive Particle Filter Algorithm at $(k+1)^{th}$ observation}
\label{alg:particleFilter}
\begin{algorithmic}
\FOR{$n = 1 \to R$}
\FOR{$\ell=0\to m-1$}
\STATE Generate $(\widetilde\Theta_{mk+\ell+1}^{(n)},\widetilde X_{mk+\ell+1}^{\epsilon,(n)})$ given $(\widetilde\Theta_{mk+\ell}^{(n)},\widetilde X_{mk+\ell}^{\epsilon,(n)})$
\ENDFOR
\STATE Compute unnormalized weight given $y_{k+1}-y_k$, $(\widetilde\Theta_{m(k+1)}^{(n)},\widetilde X_{m(k+1)}^{\epsilon,(n)})$, and $\omega_k^{(n)}$
\[\omega_{k+1}^{(n)}=\exp\left\{-\frac{1}{2}\left(\frac{y_{k+1}-y_k-\widetilde \Delta t~h\cdot\sum_{\ell=km+1}^{m(k+1)} \widetilde X_\ell^{\epsilon,(n)}}{\sqrt{\Delta t}}\right)^2\right\}\cdot\omega_k^{(n)}\]
\ENDFOR
\STATE set $c_{k+1}=\sum_n\omega_{k+1}^{(n)}$
\STATE set $\omega_{k+1}=\omega_{k+1}/c_{k+1}$

\IF{$1/\sum_n(\omega_{k+1}^{(n)})^2\leq \eta R$}
\STATE replace $\left\{(\widetilde\Theta_{m(k+1)}^{(n)},\widetilde X_{m(k+1)}^{\epsilon,(n)})\right\}_{n=1}^R$ with an $\{\omega_{k+1}^{(n)}\}_{n=1}^R$-weighted bootstrap sample of size $R$
\STATE replace $\{\omega_{k+1}^{(n)}\}_{n=1}^R$ with $(1,1,\dots,1)/R$
\ENDIF
\end{algorithmic}
\end{algorithm}

Convergence of the particle filter to the optimal filter as $R\rightarrow \infty$ can be slow, particularly if $X$ is a multi-dimensional OU process. Rao-Blackwellization can speed things up dramatically if a large portion of the hidden state's dimensionality is conditionally linear, and if the function $h$ is linear. However, if $h$ is nonlinear but $\epsilon$ is small, then the averaged-filter will be close to optimal, will be faster to compute, and its particle filter will converge faster.

In fast mean-reverting regimes, we take weights $\bar\omega$ given by
\[\bar\omega_{k+1}^{(n)}=\frac{1}{\bar c_{k+1}}\exp\left\{-\frac{1}{2}\left(\frac{y_{k+1}-y_k-\widetilde \Delta t~h\cdot\sum_{\ell=km+1}^{m(k+1)} \widetilde \Theta_\ell^{(n)}}{\sqrt{\Delta t}}\right)^2\right\}\cdot\bar\omega_k^{(n)}\]
where $\bar c_{k+1}$ is a normalizing constant. The averaged filter is
\[\mathbb E\left[g(\Theta(t_k))\Big|Y_{1:k}^0=y_{1:k}\right]\approx \sum_{n=1}^Rg(\widetilde\Theta_{km}^{(n)})\cdot\bar\omega_k^{(n)}\qquad\hbox{for }k=0,1,2,\dots,N.\]
The filter computed with the $\bar\omega$'s is appealing because it does not require samples of $\widetilde X^\epsilon$. Thus, fewer samples are required and the filter will converge at a faster rate with $R$.

\subsection{Filtering Error as a Function of $\mathbf \epsilon$}
\label{sec:filterError}
\begin{figure}[htbp] 
   \centering
   \includegraphics[scale=.58]{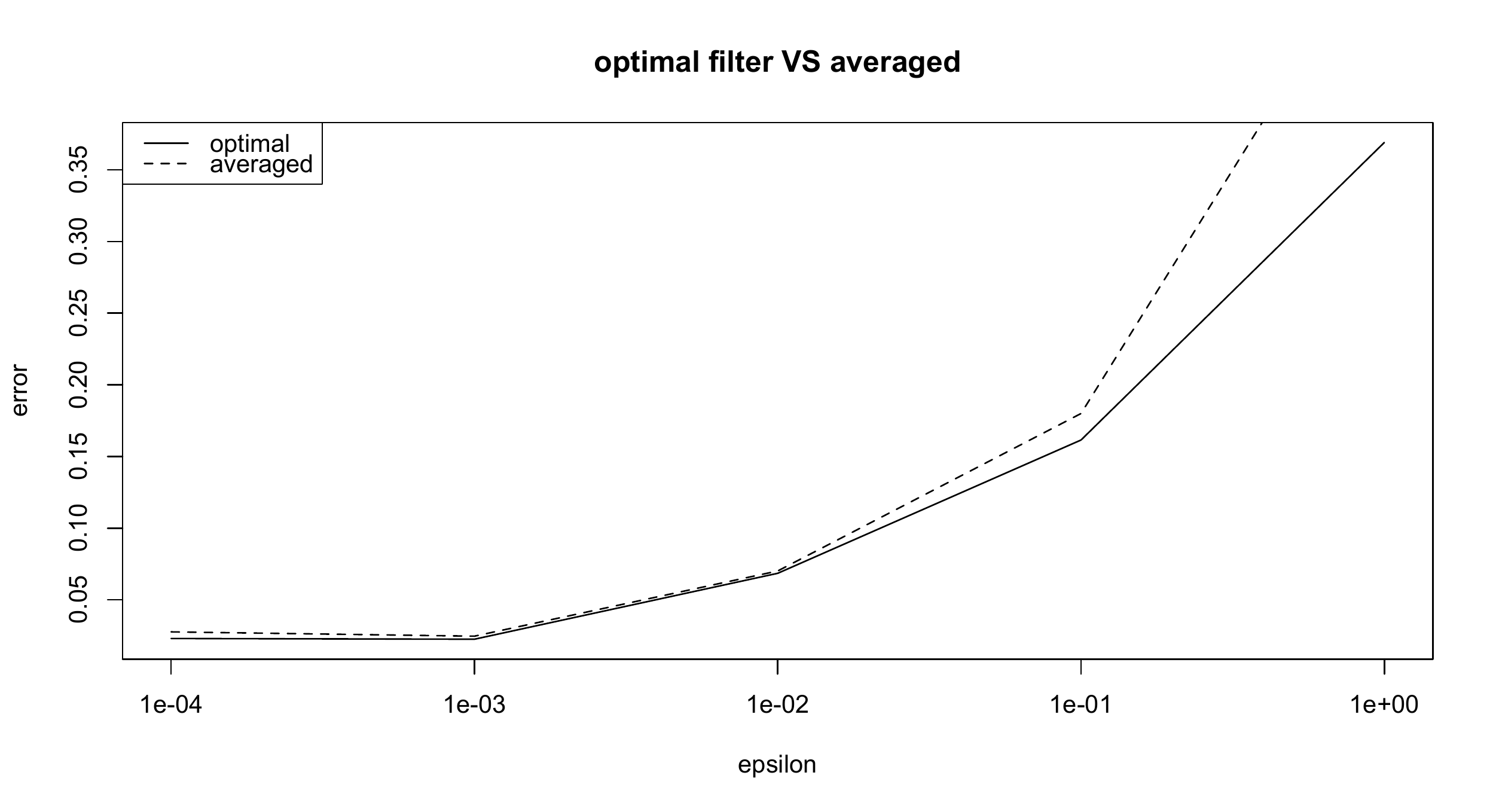} 
   \caption{\small The error of the averaged filter $\widehat\Theta_k^{avg}$ becomes closer to that of the optimal filter $\widehat\Theta_k$ as $\epsilon\searrow 0$. Also, the error of the optimal filter decreases as $\epsilon\searrow0$, indicating that SNR is increasing.}
   \label{fig:0-1error}
\end{figure}

For any $\epsilon>0$, we define $\mathbb E^\epsilon[\cdot]$ to be the expectation operator associated with the $\mathbb P^\epsilon$-measure. Then, we consider two estimates of $\Theta(t_k)$ among all the $\mathcal F_k^{Y^\epsilon}$-measurable functions: an optimal and an averaged,

\[\widehat\Theta_k^{opt} = \arg\max_{s_i}\pi_k^\epsilon(s_i)=\arg\min_{ \theta\in\mathcal F_k^{Y^\epsilon}}\mathbb E^\epsilon\mathbf 1_{\Theta(t_k)\neq \theta}\]

\[\widehat\Theta_k^{avg} = \arg\max_{s_i}\bar\pi_k(s_i).\]
For a given estimator $\widehat\Theta_k^{\{\cdot\}}$, the 0-1 error is defined and approximated as
\[error_\epsilon=\mathbb E^\epsilon \mathbf 1_{\{\Theta(t_k)\neq \widehat \Theta_k^{\{\cdot\}}\}}\approx \frac{1}{N}\sum_{k=1}^N\mathbf 1_{\{\Theta(t_k)\neq \widehat \Theta_k^{\{\cdot\}}\}}\qquad\hbox{for $N$ large.}\]
From the theory presented in the Sections leading up to now, we should find that the 0-1 error of the averaged estimate converges to the error of the optimal estimate as $\epsilon\searrow 0$. The estimators are computed with the following approximated $\pi$'s:
\[\pi_k^\epsilon(s_i) \approx \sum_{n=1}^R\omega_k^{(n)}\cdot\mathbf 1_{\{\tilde\Theta_{km}^{(n)}=s_i\}}\]
\begin{center}and\end{center}
\[\bar\pi_k(s_i) \approx \sum_{n=1}^R\bar\omega_k^{(n)}\cdot\mathbf 1_{\{\tilde\Theta_{km}^{(n)}=s_i\}}.\]
The figures in this section and the next are generated with the following parameterization of the problem:
\[\Delta t=.01,\qquad N = 100,000\qquad\mathcal S = [-3.\overline{333},3.\overline{333}],\qquad\alpha=10,\qquad\beta=5\]
\[Q=\left(\begin{array}{cc}
-\alpha&\alpha\\
\beta&-\beta
\end{array}\right)\qquad X(0) \sim U[-1,1],\qquad  h(x)=10\cdot x.\]
For $R=100$ and $m=5$ (so that $\widetilde \Delta t=.0005$), Figure \ref{fig:0-1error} shows how the 0-1 error of the average and optimal converge. The figure was computed using the Rao-Blackwellized filter described in Appendix \ref{app:raoBlackwell}, rather than the particle filter of Algorithm \ref{alg:particleFilter}.  From the figure, we can clearly see that there is convergence of the averaged filter's error to the optimal filter's error as a $\epsilon$ decreases,
\[0\leq  \mathbb E^\epsilon\left(\mathbf 1_{\{\Theta_{t_k}\neq\widehat\Theta_k^{avg}\}}-\mathbf 1_{\{\Theta_{t_k}\neq\widehat\Theta_k^{opt}\}}\right)\stackrel{\epsilon\searrow 0}{\longrightarrow}0.\]
It should also be noted that we can easily identify the signal-to-noise ratio (SNR) for this problem because it is linear. In fact, in Figure \ref{fig:0-1error} we see a decrease in the optimal filter's error as $\epsilon$ decreases, indicating that the variance of the estimator decreases as well. 

\textcolor{black}{It is interesting to compare similarities and differences between the filter in \cite{papanicolaou2010} (described in Section \ref{sec:papanicolaou2010}) and the one in this paper. The two filters are quite different because they have different scalings, and because the filter in \cite{papanicolaou2010} is derived from a strong-sense limit whereas the one in this paper is derived under the weak-topology. Also, averaging over $h(x)$ makes the filter in this paper simpler and easier to implement than the one in \cite{papanicolaou2010} where no such averaging occurs. However, these filters are similar because both of their performances are contingent on the goodness of the asymptotic approximation. Indeed, the ability of the reduced-dimension filter to improve as the rate of mean-reversion increases is seen in the simulations of both papers. In fact, \cite{papanicolaou2010} has a plot that is similar to Figure \ref{fig:0-1error}.}

\textcolor{black}{Our conclusion from this section is that the averaged-filter can perform well, but its closeness to the optimal filter will depend on $\epsilon$. Whether or not optimality is essential will depend on the problem and the user's specific needs; it might be merely an extra 1\%-2\% in error that is under scrutiny.}

\subsection{Search for the Optimal $\mathbf \Delta t$}
\label{sec:empirical}

\begin{figure}[ht]

\centering
\includegraphics[scale=.48]{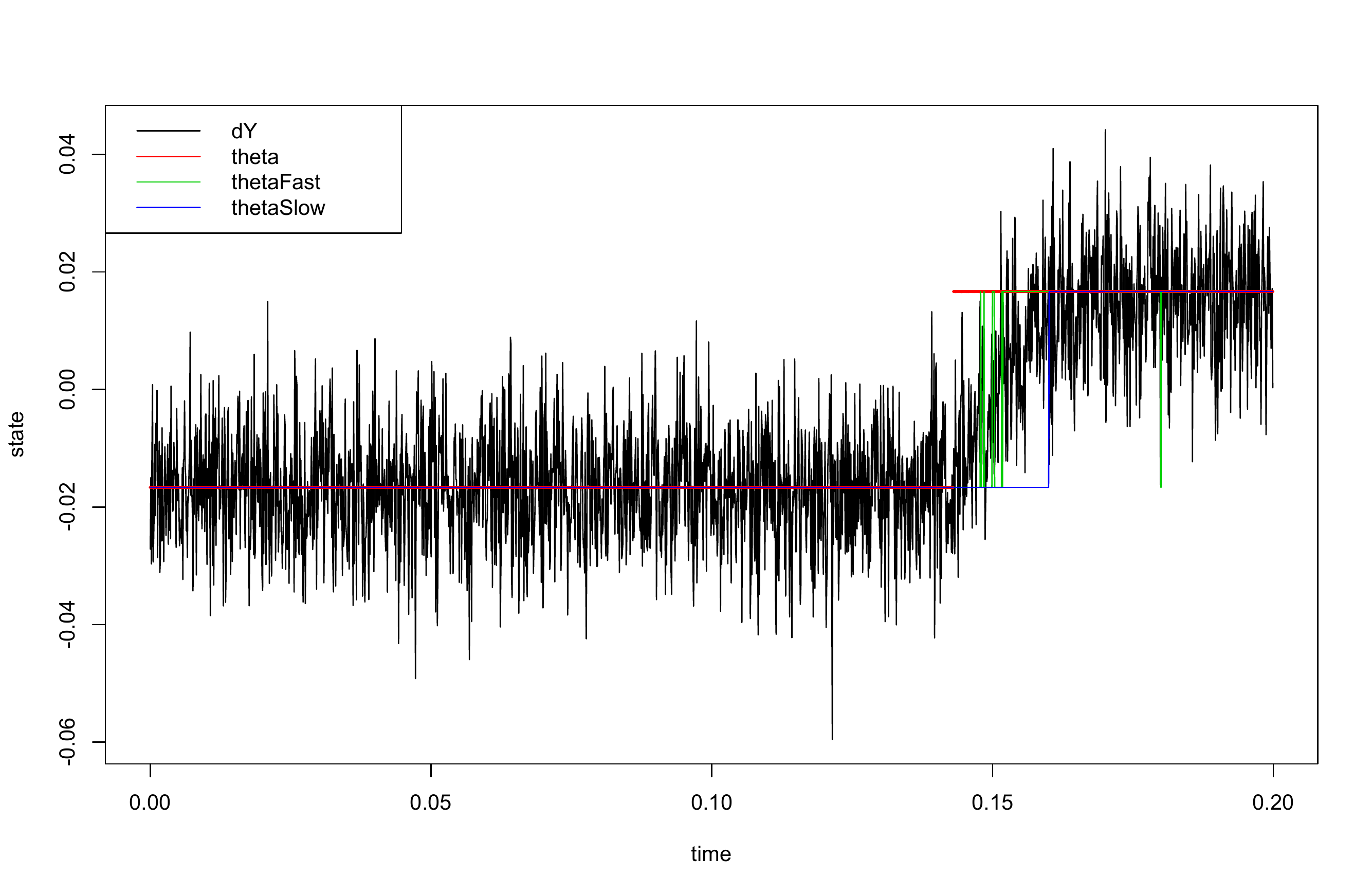}

\caption{\small A realization where sampling too slow $\Delta t=.01$ or too fast $\Delta t=.0001$ will not be optimal. The black line is the observations $\Delta Y_k$, the red line is the true state $\Theta$, the green line is the filter obtain by sampling observations very slow, and the blue line is the filter obtained by sampling too fast. Notice how the green line sometimes tracks noise.}
\label{fig:trackingNoise}

\end{figure}

\begin{figure}[ht]
\centering
\includegraphics[scale=.58]{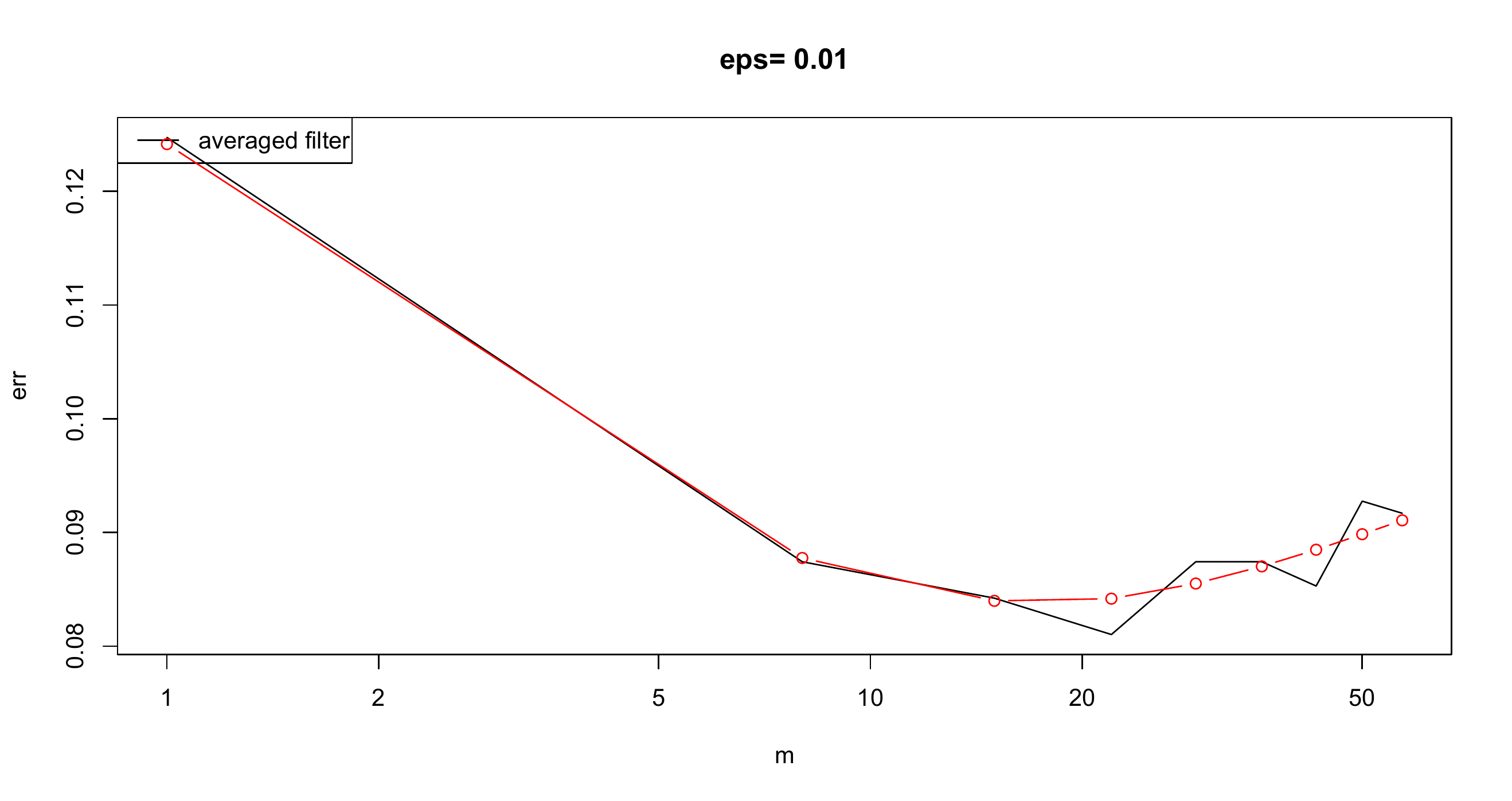}
\caption{\small The 0-1 error for $\epsilon=.01$ and varying $\Delta t=m\tilde\Delta t=m\times 10^{-4}$. The solid line is the error from the simulation, and the dotted line is a spline included so that the low-point is clearer.}
\label{fig:empiricalDelta}
\end{figure}
In this section we run simulations in an attempt to find the optimal rate at which observations should be sampled. The simulation and parameters are the same as they were in Section \ref{sec:filterError}, except we will vary $\Delta t$. The emphasis will be on the choice of $\Delta t$ which minimizes filtering error,

\begin{equation}\label{eq:optimizationCriterion}
\min_{\Delta t}\left(\mathbb E^\epsilon\mathbf 1_{\Theta(t_k)\neq \widehat\Theta_k^{avg}}\right).
\end{equation}
Simply stated, sampling too fast will eliminate any hypothesis regarding a fast mean-reverting regime, while sampling too slow will cause the filter's memory to deteriorate and the performance will be poor. Therefore, given all the other parameters, the user must choose the value of $\Delta t$ that gives the averaged filter a good chance to perform well.

The first thing to notice is the potential for the filter to track noise. In Figure \ref{fig:trackingNoise}, two filters are plotted alongside the true state and the observations. One of the filters takes $\Delta t=.01$ and returns an estimate that misses the transition in $\Theta$. The other filter takes $\Delta t=.0001$, which is entirely too fast for the fast-mean-reverting hypothesis and results in the tracking of noise.

Ideally, the averaged filter will be implemented with a priori knowledge of the optimal sampling rate, but analytical expressions for such optima are not known. Instead, the simulation is run several times with different $\Delta t$'s, and the errors are compared. Figure \ref{fig:empiricalDelta} shows the optimal $\Delta t$ to be around $.002$ for this model ($\Delta t=m\tilde\Delta t=m\times .0001=.002$). For each $m$ in Figure \ref{fig:empiricalDelta}, the averaged filter is computed with $R= 20\times\log_2(m)$ as the number of particles.

The conclusion is that $\Delta t$ cannot be too small; otherwise, the fast filter becomes more erroneous. Intuitively, it makes sense if one understands that some `coarse-graining' is required to insure that the fast-averaging of the model's integrals will occur between observation times. For this reason, the average filter does not apply for filters based on a continuum of observations. In practice, the user must consider the model's parameters and make a decision on $\Delta t$'s value.
\section{Conclusion} We have considered a filtering problem for a hidden Markov model where one of the hidden states is fast mean-reverting with a shifting mean given by a Markov chain that is also a hidden process that models regime change. We have derived an averaged filter of low complexity for tracking directly the regime changes. We have given examples of where this model may be useful and we have also provided some related examples to show that subtle changes to the model will change the asymptotic behavior of the filter. We have also given a rigorous treatment and proof to show that the averaged filter is asymptotically optimal. 

We have also presented the results of some numerical simulations to show how the averaged filter performs as the rate of mean reversion increases. These simulations indicate that the averaged filter is a powerful numerical method if used correctly. Specifically, the filter's performance depends in large part on parameter values for the rate of mean-reversion and the rate of observations.

In the future, possible research on nonlinear filters for HMMs with fast mean-reverting states includes robustness analysis, parameter estimation and the EM algorithm, and further application of the algorithms to financial data.

\appendix

\section{Tightness of Marginal Measures and Weak Convergence to a Unique Limit}
\label{sec:details}
This section contains the lemmas and theorems necessary to prove (\ref{eq:weakLimit}). In doing so, it is useful to introduce a marginal probability space $(\bar\Omega,\bar{\mathcal F}_t,\bar{\mathbb P})$ for the non-Markov process $(\Theta^\epsilon(t),V^\epsilon(t))$, where
 \[\bar\Omega = D\left([0,T];\mathcal S\right)\times C\left([0,T];\mathbb R\right)\]
and we define the measures 

\[\bar{\mathbb P}^\epsilon(A) = \bar{\mathbb P}((\Theta^\epsilon(\cdot),V^\epsilon(\cdot))\in A)\]
 for all $\epsilon>0$ and $A\in \bar{\mathcal F}_T$. It can be shown that the process $\Theta^\epsilon(\cdot)$ is tight  in the space $D([0,T];\mathcal S)$ and $V^\epsilon(\cdot)$ is tight in the space $C([0,T];\mathbb R)$, which therefore makes the family $(\bar{\mathbb P}^\epsilon)_{\epsilon>0}$ tight in $\bar\Omega$. Furthermore, for any bounded and differentiable $g:\mathcal S\times\mathbb R\rightarrow\mathbb R$, we can define the following continuous functional for any path $(\theta(\cdot),v(\cdot))\in\bar\Omega$,

\begin{equation}
\label{eq:martingale}
 \phi_g(t)  = g\left(\theta(t),v(t)\right)-g\left(\theta(0),v(0)\right)-\int_0^t\left(\bar Q+\bar h_{\theta(\tau)}\frac{\partial}{\partial v}\right)g(\theta(\tau),v(\tau))d\tau
\end{equation}
Under the measure $\bar{\mathbb P}^\epsilon$ the expectation of $M_g(t)$ is defined as

\[\mathbb E^{\bar{\mathbb P}^\epsilon}[\phi_g(t)] \doteq \int \phi_g(t)d\mathbb P^\epsilon\]
for which we will show that for any $t,s\in [0,T]$ with $t\geq s$,
\[\mathbb E^{\bar{\mathbb P}^\epsilon}\left[\phi_g(t)-\phi_g(s)\Big|\mathcal F_s\right]\rightarrow0\qquad\hbox{as }\epsilon\searrow0\]
in probability, which is sufficient to conclude that for any weakly convergent subsequence of $(\bar{\mathbb P}^{\epsilon})_{\epsilon>0}$ with limit $\bar{\mathbb P}$, we have that $\phi_g$ is a $\bar{\mathbb P}$-Martingale. 
Therefore, from the tightness of the marginal measures $(\bar{\mathbb P}^\epsilon)_{\epsilon>0}$ it follows that
\[\left(\Theta^\epsilon(\cdot),V^\epsilon(\cdot)\right)\Rightarrow (\overline\Theta(\cdot),\overline V(\cdot))\qquad\hbox{as }\epsilon\searrow0\]
where $ (\overline\Theta(\cdot),\overline V(\cdot))$ uniquely solves the $\bar{\mathbb P}$-Martingale problem associated with $\phi_g$ in (\ref{eq:martingale}), and is a Markov process with generator $\bar Q+\bar h_{s_i}\frac{\partial}{\partial v}$. 
 
To summarize, the goal of this section is to prove (\ref{eq:weakLimit}) by showing that the marginal measures $(\bar{\mathbb P}^\epsilon)_{\epsilon>0}$ are tight in $\bar\Omega$ and that every convergent subsequence must converge to the measure of a process $(\overline\Theta(\cdot),\overline V(\cdot))$ which is the unique solution to a specific Martingale problem. This is a well known method of proving weak convergence and is presented in \cite{yinZhangCont}, for example.
\subsection{Tightness} General result regarding tightness of measures and weak convergence can be found in the book by Ethier and Kurtz \cite{ethierKurtz}, (chapter 3). In fact the proof of the tightness of $\Theta^\epsilon(\cdot)$ that is shown in this section is a specific case of lemma 6.1 on page 122 of their book.

 Tightness of a family of measures means that for all $\delta >0$ there exists a compact set $K_\delta \subset \bar\Omega$ such that
\[\sup_\epsilon\bar{\mathbb P}^\epsilon(K_\delta^c)<\delta\]
It turns out that because of the bound in (\ref{eq:upperJumpRateBound}), the family of measures is tight in $\bar\Omega$.

\begin{lemma}\label{lemma:tightness} Assuming that the initial condition $V_0$ is almost-surely bounded, the bound in (\ref{eq:upperJumpRateBound}) insures that $(\bar{\mathbb P}^\epsilon)_{\epsilon>0}$ is tight in $\bar\Omega$.
\end{lemma}

\textbf{Proof of Lemma \ref{lemma:tightness}:} First, let's prove that $V^\epsilon(\cdot)$ is tight in $C([0,T];\mathbb R)$. Since there exists a positive constant $\kappa<\infty$ such that $\mathbb P(|V_0|<\kappa)=1$ and because $h$ is bounded, the function $V^\epsilon(t)$ is uniformly bounded,
\[|V^\epsilon(t)|\leq \kappa+T\|h\|_\infty\qquad\forall t\in[0,T],\forall \epsilon>0\]
Furthermore, $V^\epsilon(t)$ is uniformly equicontinuous with modulus of continuity $\|h\|_\infty$,
\[|V^\epsilon(t)-V^\epsilon(s)|\leq \|h\|_\infty|t-s|\qquad\forall t,s\in[0,T],\forall \epsilon>0\]
Therefore, by the Ascoli theorem, any subsequence of $V^\epsilon(\cdot)$ has a further subsequence that converges uniformly, and therefore the support of all measures $\bar{\mathbb P}^\epsilon$ on $V^\epsilon$ is compact. Therefore, for any $\delta >0$ there is a set $K_\delta(V)$ such that 
\[\bar{\mathbb P}^\epsilon(K_\delta(V))>1-\delta/2\qquad\forall\epsilon>0.\]

Next, to prove compactness of $\Theta^\epsilon(\cdot)$, consider any $f\in D([0,T],\mathcal S)$ and define the switching times $\{\tau_k(f)\}_{k=0}^\infty$ for $k=1,2,3,...$
\[\tau_k(f) = \bigg\{\begin{array}{cl}
\inf\{s\in(\tau_{k-1}(f),T]:f(s)\neq f(s^-)\hbox{ or }s=T\}&\hbox{if }\tau_{k-1}(f)<T\\
T&\hbox{otherwise}
\end{array}\]
with $\tau_0(f)=0$. For any compact set $U\subset \mathbb R$, let the set $A(U,\delta)$ be defined as
\[A(U,\delta) = \{f:f(s)\in U\hbox{ for all }s\leq T,\tau_k(f)-\tau_{k-1}(f)>\delta,\forall k\hbox{ with }\tau_k(f)<T\}\]
The closure of $A(U,\delta)$ is compact:\\

Given any sequence $f_n\in A(U,\delta)$, by a diagonal argument there exists a subsequence $f_{n_\ell}$, and limits $t_k$ and $a_k$ such that
\[\tau_k(f_{n_\ell})\rightarrow t_k\in[0,T]\]
and 
\[f_{n_\ell}(\tau_k(f_{n_\ell}))\rightarrow a_k\in U\]
as $\ell\nearrow\infty$. Therefore, 
\[f_{n_\ell}(s)\rightarrow a_k,\qquad\hbox{for }s\in[t_k,t_{k+1})\]
which proves that $\overline{A(U,\delta)}$ is sequentially compact, and since $D([0,t];\mathcal S)$ is equipped with the point-wise metric, the subset is therefore compact.\\

For the Markov process $\Theta^\epsilon(\cdot)$, take $U$ s.t. $\{s_1,....,s_M\}\subset U$ and recall that
\[\mathbb P^\epsilon(A(U,\delta)^c)\leq 1-e^{-\beta\delta}\leq \beta\delta\qquad\forall \epsilon>0\]
which shows that $\Theta^\epsilon(\cdot)$ is compact in $D([0,T];\mathbb R)$.

Finally, tightness of $(\bar{\mathbb P}^\epsilon)_{\epsilon>0}$ follows by considering the compact set $K_\delta(V)\times\overline {A(U,\delta/2\beta)}$,
\[\bar{\mathbb P}^\epsilon\left((K_\delta(V)\times\overline {A(U,\delta/2\beta)})^c\right)\leq\bar{\mathbb P}^\epsilon\left(K_\delta(V)^c\right)+\bar{\mathbb P}^\epsilon\left(\overline {A(U,\delta/2\beta)}^c\right) <\delta/2+\delta/2=\delta\]

$\blacksquare$\\
\subsection{Convergence of the Martingale Problem} 
\label{app:convergenceMartingale}

In this section it will be shown by theorem \ref{thm:meanConvergence} that for any $t,s\in[0,T]$ with $t\geq s$, 
\[\mathbb E^{\mathbb P^\epsilon}\left[\phi_g(t)-\phi_g(s)|\mathcal F_s\right]\rightarrow 0\qquad\hbox{ as }\epsilon\searrow 0\]
in probability. 


To prove that the family of marginal measures will converge to the unique solution of the Martingale problem, the following lemma will be necessary:

\begin{lemma}\label{lemma:errorBound}
For any positive $\Delta t>0$ and any function $g$ such that $\left\|g\right\|_\infty+\left\|\frac{\partial}{\partial v}g\right\|_\infty\leq C$, the expectation of the following terms converge to a term of $o(\Delta t)$ as $\epsilon\searrow 0$
\[\lim_{\epsilon\searrow0}\left|\mathbb E\left[\int_0^{\Delta t}\left(Q(X^\epsilon(s))-\bar Q\right)g(\Theta^\epsilon(s),V^\epsilon(s))ds\Big|\mathcal F_0\right]\right|\leq2\beta^2 \Delta t^2\]
\[\lim_{\epsilon\searrow0}\left|\mathbb E\left[\int_0^{\Delta t}\left(h(X^\epsilon(s))-\bar h_{\Theta^\epsilon(s)}\right)\frac{\partial}{\partial v}g(\Theta^\epsilon(s),V^\epsilon(s))ds\Big|\mathcal F_0\right]\right|\leq2\beta^2\Delta t^2\]
in probability.\\
\end{lemma}
\noindent\textbf{Proof of lemma \ref{lemma:errorBound}:} Start by proving the first statement. By iteratively conditioning the $\sigma$-algebras, the expectation can be separated into two parts,
\[\left|\mathbb E\left[\int_0^{\Delta t}\left(Q(X^\epsilon(s))-\bar Q\right)g(\Theta^\epsilon(s),V^\epsilon(s))ds\Big|\mathcal F_0\right]\right|\]

\[\leq\left|\mathbb E\left[\int_0^{\Delta t}\left(Q(X^\epsilon(s))-\bar Q\right)g(\Theta^\epsilon(s),V^\epsilon(s))ds\Bigg|\mathcal F_0\cap\{\Theta^\epsilon(s) \equiv\Theta^\epsilon(0)~~\forall s\in[0,\Delta t]\}\right]\right|\]
\[\qquad\qquad\qquad\qquad\times\mathbb P(\hbox{$\Theta^\epsilon(s)\equiv\Theta^\epsilon(0)$ for $s\in[0,\Delta t]$}) \]

\[+\left|\mathbb E\left[\int_0^{\Delta t}\left(Q(X^\epsilon(s))-\bar Q\right)g(\Theta^\epsilon(s),V^\epsilon(s))ds\Bigg|\mathcal F_0\cap\{\max_{s\leq\Delta t}|\Theta^\epsilon(s)-\Theta^\epsilon(0)|>0\}\right]\right|\]

\[\qquad\qquad\qquad\qquad\qquad\times\left(1-\mathbb P(\hbox{$\Theta^\epsilon(s)\equiv\Theta^\epsilon(0)$ for $s\in[0,\Delta t]$})\right)\]
Now, it obvious that $\mathbb P(\hbox{$\Theta^\epsilon(s)\equiv\Theta^\epsilon(0)$ for $s\in[0,\Delta t]$}) \leq 1$, and from the bound in (\ref{eq:upperJumpRateBound}) it follows that

\[1-\mathbb P(\hbox{$\Theta^\epsilon(s)\equiv\Theta^\epsilon(0)$ for $s\in[0,\Delta t]$}) \leq 1-e^{-\beta\Delta t}\]
On intervals where $\Theta^\epsilon(t)$ is constant, we have
\[X^\epsilon(t)=_d \Theta^\epsilon(0)+e^{-t/\epsilon}(x-\Theta^\epsilon(0))+N\left(0,\frac{1}{2}\left(1-e^{-2t/\epsilon}\right)\right)\]
and on these intervals the densities of  $X^\epsilon(t)$ given $X(0)=x$ can be written as $\mu_{t/\epsilon}(y|x)\in\mathbb R^{M\times M}$ as
\[\mu_{t/\epsilon}(y|x)=\frac{1}{\sqrt{\pi(1-e^{-2t/\epsilon})}}\left[\begin{array}{ccc}
\exp\{-\frac{(y-s_1-e^{-2t/\epsilon}(x-s_1))^2}{1-e^{-2t/\epsilon}}\}&\dots&0\\
&\ddots&\\
0&\dots&\exp\{-\frac{(y-s_M-e^{-2t/\epsilon}(x-s_M))^2}{1-e^{-2t/\epsilon}}\}
\end{array}\right]
\]
From here, there is the following upper-bound,
\[\left|\mathbb E\left[\int_0^{\Delta t}\left(Q(X^\epsilon(s))-\bar Q\right)g(\Theta^\epsilon(s),V^\epsilon(s))ds\Big|\mathcal F_0\right]\right|\]

\[\leq\left|\mathbb E\left[\int_0^{\Delta t}\left(Q(X^\epsilon(s))-\bar Q\right)g(\Theta^\epsilon(s),V^\epsilon(s))ds\Bigg|\mathcal F_0\cap\{\Theta^\epsilon(s) \equiv\Theta^\epsilon(0)~~\forall s\in[0,\Delta t]\}\right]\right|\]
\[\qquad\qquad\qquad\qquad+2\beta\Delta t(1-e^{-\beta\Delta t})\]

\[=\left|\int_0^{\Delta t}\left(\int \mu_{s/\epsilon}(y|X^\epsilon(0))Q(y)dy-\bar Q\right)g(\Theta^\epsilon(0),V^\epsilon(0))ds\right|+2\beta\Delta t(1-e^{-\beta\Delta t})\]

\[\leq\|g\|_\infty\left\|\int_{\epsilon^\alpha}^{\Delta t}\underbrace{\left(\int \mu_{s/\epsilon}(y|X^\epsilon(0))Q(y)dy-\bar Q\right)}_{:=\gamma(\Delta t,\epsilon)}ds\right\|+\underbrace{2\beta\Delta t(1-e^{-\beta\Delta t})}_{\leq 2\beta^2\Delta t^2}+\underbrace{\|g\|_\infty\sup_x\|Q(x)-\bar Q\|}_{=C'}\epsilon^\alpha\]
where $0<\alpha<1$. It is clear the $\gamma(\Delta t,\epsilon)\rightarrow 0$ as $\epsilon\searrow 0$ for all $\Delta t>0$ in probability. Therefore, we have

\[\left|\mathbb E\left[\int_0^{\Delta t}\left(Q(X^\epsilon(s))-\bar Q\right)g(\Theta^\epsilon(s),V^\epsilon(s))ds\Big|\mathcal F_0\right]\right|\leq\Delta t\gamma(\Delta t,\epsilon) +2\beta^2\Delta t^2+C'\epsilon^\alpha\rightarrow2\beta^2 \Delta t^2\]
which proves the first statement of the lemma.

Repeating these steps with $Q(X^\epsilon(s))$ replaced by $h(X^\epsilon(s))$ and $\bar Q$ replace by $\bar h_{\Theta^\epsilon(s)}$, the second statement in the lemma is proven provided that there is a constant such that 
\[\left\|g\right\|_\infty+\left\|\frac{\partial}{\partial v}g\right\|_\infty\leq C.\]

$\blacksquare$\\
With lemma \ref{lemma:errorBound}, the key calculation for proving weak convergence of $(\Theta^\epsilon(\cdot),V^\epsilon(\cdot))$ can be shown. Namely, from the time-homogeneity and Markov property of the solution to the Martingale problem associated with $\phi_g$, it is sufficient to show that the expectation under $\mathbb P^\epsilon$ of the function $\phi_g(t)$ conditioned on $\mathcal F_0$ converges to zeros as $\epsilon\searrow 0$.
\begin{theorem}\label{thm:meanConvergence}
For any function $g(\theta,v)$ with $\left\|g\right\|_\infty+\left\|\frac{\partial}{\partial v}g\right\|_\infty\leq C$, condition (\ref{eq:upperJumpRateBound}) and the boundedness of $h$ insure that for any $t,s\in[0,T]$ with $t\geq s$, 
\[\mathbb E^{\mathbb P^\epsilon}[\phi_g(t)-\phi_g(s)|\mathcal F_s]\rightarrow 0\]
as $\epsilon\searrow 0$ in probability.

\end{theorem}
\noindent\textbf{Proof of Theorem \ref{thm:meanConvergence}:} Without loss of generality we take $s=0$ and any $t\in[0,T]$. For a fixed step-size $\Delta t>0$, we consider the following discretization of the interval $[0,t]$,
\[0<\Delta t<2\Delta t<\dots \dots<N\Delta t = t\]
For any $n<N$, let $t_n = n\Delta t$. We then have a collapsing sum
\[\mathbb E^{\mathbb P^\epsilon}[\phi_g(t)|\mathcal F_0]\]

\[=\mathbb E\left[g(\Theta^\epsilon(t),V^\epsilon(t))-g(\Theta^\epsilon(0),V^\epsilon(0))|\mathcal F_0\right]\qquad\qquad\qquad\qquad\qquad\qquad\qquad\qquad\qquad\qquad\]

\[\qquad\qquad-\mathbb E\left[\int_0^t\left(\bar Q+\bar h_{\Theta^\epsilon(s)}\frac{\partial}{\partial v}\right)g(\Theta^\epsilon(s),V^\epsilon(s))ds\Bigg|\mathcal F_0\right]\]

\[=\mathbb E\left[\sum_{n=0}^{N-1}g(\Theta^\epsilon(t_{n+1}),V^\epsilon(t_{n+1}))-g(\Theta^\epsilon(t_n),V^\epsilon(t_n))|\mathcal F_0\right]\qquad\qquad\qquad\qquad\qquad\qquad\qquad\qquad\qquad\]

\[\qquad\qquad-\mathbb E\left[\sum_{n=0}^{N-1}\int_{t_n}^{t_{n+1}}\left(\bar Q+\bar h_{\Theta^\epsilon(s)}\frac{\partial}{\partial v}\right)g(\Theta^\epsilon(s),V^\epsilon(s))ds\Bigg|\mathcal F_0\right]\]

\[=\underbrace{\mathbb E\left[\sum_{n=0}^{N-1}\mathbb E[g(\Theta^\epsilon(t_{n+1}),V^\epsilon(t_{n+1}))|\mathcal F_{t_n}]-g(\Theta^\epsilon(t_n),V^\epsilon(t_n))\Big|\mathcal F_0\right]}_{(*)}\qquad\qquad\qquad\qquad\qquad\qquad\qquad\qquad\qquad\qquad\qquad\qquad\]

\[-\underbrace{\mathbb E\left[\sum_{n=0}^{N-1}\int_{t_n}^{t_{n+1}}\mathbb E\left[\left(\bar Q+\bar h_{\Theta^\epsilon(s)}\frac{\partial}{\partial v}\right)g(\Theta^\epsilon(s),V^\epsilon(s))\Big|\mathcal F_{t_n}\right]ds\Bigg|\mathcal F_0\right]}_{(**)}\]
The summand in $(**)$ can be shown to be equal to the summand in $(*)$ plus a correction term,

\[\int_{ t_n}^{t_{n+1}}\mathbb E \left[\left(\bar Q+\bar h_{\Theta^\epsilon(s)}\frac{\partial}{\partial v}\right)g(\Theta^\epsilon(s),V^\epsilon(s))\Big|\mathcal F_{t_n}\right]ds \]

\[=\underbrace{\int_{ t_n}^{t_{n+1}}\mathbb E \left[\left(\bar Q+\bar h_{\Theta^\epsilon(s)}\frac{\partial}{\partial v}-Q(X^\epsilon(s))-h(X^\epsilon(s))\frac{\partial}{\partial v}\right)g(\Theta^\epsilon(s),V^\epsilon(s))\Big|\mathcal F_{t_n}\right]ds}_{=e_n(\epsilon,\Delta t))} \]

\[+\int_{ t_n}^{t_{n+1}}\mathbb E \left[\left(Q(X^\epsilon(s))-h(X^\epsilon(s))\frac{\partial}{\partial v}\right)g(\Theta^\epsilon(s),V^\epsilon(s))\Big|\mathcal F_{t_n}\right]ds\]

\[=\mathbb E[g(\Theta^\epsilon(t_{n+1}),V^\epsilon(t_{n+1}))|\mathcal F_{t_n}]-g(\Theta^\epsilon(t_n),V^\epsilon(t_n))+e_n(\epsilon,\Delta t)\]
By lemma \ref{lemma:errorBound}, the correction term $e_n(\epsilon,\Delta t)$ can be controlled so that the limit as $\epsilon\searrow 0$ is $o(\Delta t)$ in probability. Therefore, 

\[\left|\mathbb E^{\mathbb P^\epsilon}[\phi_g(t)|\mathcal F_0]\right|\leq \sum_{n=0}^{N-1}|e_n(\epsilon,\Delta t)|\rightarrow\sum_{n=0}^{N-1}o(\Delta t) = O(\Delta t)\qquad\hbox{as }\epsilon\searrow 0\]
and the entire sum converges to something of order $O(\Delta t)$ as $\epsilon\searrow 0$ in probability, which proves the theorem since $\Delta t$ is arbitrarily small.\\
$\blacksquare$ 

Theorem \ref{thm:meanConvergence} confirms that any convergent sequence in the family of tight measure $(\bar{\mathbb P}^\epsilon)_{\epsilon>0}$ must have a limit $\bar{\mathbb P}$ which is the measure for which paths in the space $\bar\Omega$ uniquely solve the Martingale problem associated with $\phi_g$ . Therefore, this proves (\ref{eq:weakLimit}) and the reduced filter of theorem \ref{thm:averagedFilter} stands.
\section{The Rao-Blackwellized Filter for Linear Observations, $(\mathbf{h(x)=h\cdot x})$} 
\label{app:raoBlackwell}
The Rao-Blackwellized filter for a general class of target tracking problems is covered in \cite{gustafsson2002,gustafsson2005}. The Rao-Blackwellized filter will be specific to the model in question; this appendix covers the case for the OU model of this paper.

Using auto-regressive coefficient $a=e^{-\widetilde\Delta t/\epsilon}$, we have the following recursion for Algorithm \ref{alg:samples}:
\[\widetilde X_{(k+1)m}^\epsilon = a\widetilde X_{(k+1)m-1}^\epsilon+(1-a)\widetilde\Theta_{(k+1)m}^\epsilon+\sqrt{\frac{1-a^2}{2}}\mathcal W_{(k+1)m}\]

\[=a^2\widetilde X_{(k+1)m-2}^\epsilon+(1-a)\widetilde\Theta_{(k+1)m}^\epsilon+(1-a)a\widetilde\Theta_{(k+1)m-1}^\epsilon+\sqrt{\frac{1-a^2}{2}}\left(\mathcal W_{(k+1)m}+a\mathcal W_{(k+1)m-1}\right)\]

\[\dots\dots\dots\]

\[=a^m\widetilde X_{km}^\epsilon+(1-a)\sum_{\ell=1}^ma^{m-\ell} \widetilde\Theta_{km+\ell}+\sqrt{\frac{1-a^2}{2}}\sum_{\ell=1}^ma^{m-\ell} \mathcal W_{km+\ell},\]
where $\mathcal W_\ell\sim iid N(0,1)$. Using this recursive formula, the values of $\widetilde X^\epsilon$ and $\widetilde\Theta$ that occur at the $m$-many times between observations can be stored in a vectors $\vec {\tilde X}_{k+1}^\epsilon=(\widetilde X_{(k+1)m}^\epsilon,\widetilde X_{(k+1)m-1}^\epsilon,\dots,\widetilde X_{km+1}^\epsilon)^T$, and $\vec{\tilde\Theta}_{k+1}=(\widetilde \Theta_{(k+1)m},\widetilde\Theta_{(k+1)m-1},\dots,\widetilde \Theta_{km+1})^T$ where `$T$' denotes matrix/vector transpose. There matrix/vector form of the recursion is
\[\vec{\tilde X}_{k+1}^\epsilon=A\vec{\tilde X}_k^\epsilon+B\vec{\tilde\Theta}_{k+1}+R\vec {\mathcal W}_{k+1}\]
where $\vec{\mathcal W}_{k+1}=(\mathcal W_{(k+1)m},\mathcal W_{(k+1)m-1},\dots, \mathcal W_{km+1})^T$, and the matrices are

\[A=\left(\begin{array}{cccc}
a^m&0&\dots&0\\
a^{m-1}&0&\dots&0\\
\vdots&\vdots&\ddots&\vdots\\
a&0&\dots&0\\
\end{array}\right),
\qquad B=(1-a)\left(\begin{array}{ccccc}
1&a&\dots&a^{m-2}&a^{m-1}\\
0&1&\dots&a^{m-3}&a^{m-2}\\
\vdots&\vdots&\ddots&\vdots&\vdots\\
0&0&\dots&1&a\\
0&0&\dots&0&1
\end{array}\right)\]

\[\hbox{and}\qquad R=\sqrt{\frac{1-a^2}{2}}\left(\begin{array}{ccccc}
1&a&\dots&a^{m-2}&a^{m-1}\\
0&1&\dots&a^{m-3}&a^{m-2}\\
\vdots&\vdots&\ddots&\vdots&\vdots\\
0&0&\dots&1&a\\
0&0&\dots&0&1
\end{array}\right).\]
Observations are given by
\[Y_{k+1}^\epsilon-Y_k^\epsilon=H\vec{\tilde X}_{(k+1)m}^\epsilon+\Delta Z_k\]
with $\Delta Z_k\sim iid N(0,\widetilde \Delta t)$, and $H= \Delta t\cdot(h,h,h,\dots,h)$.

The idea of Rao-Blackwellization is to generate particles of $\widetilde\Theta$, and then exploit the remaining linearity. For each sample-path of $\widetilde\Theta_{1:km}^{(n)}$, there corresponds a weight $\omega_k^{(n)}$, and a Kalman filter $(\widehat{\vec{\tilde X}}_k^{\epsilon,(n)},\Sigma_k)$ where
\begin{eqnarray}
\nonumber
\widehat{\vec{\tilde X}}_k^{\epsilon,(n)}&=&\mathbb E[X_k^\epsilon|\mathcal F_k,\vec{\tilde\Theta}_{1:km}^{(n)}]\\
\nonumber
\Sigma_k&=&\mathbb E(\vec{\tilde X}_{km}^\epsilon-\widehat{\vec{\tilde X}}_k^{\epsilon,(n)})(\vec{\tilde X}_{km}^\epsilon-\widehat{\vec{\tilde X}}_k^{\epsilon,(n)})^T.
\end{eqnarray}
The covariance matrix $\Sigma_k$ evolves independently of $Y$, and it also evolves independently of $\widetilde\Theta^{(n)}$ for the particular model that were implementing. Therefore, $\Sigma_k$ does not depend on the data or the samples. However, $\widehat{\vec{\tilde X}}_k^{\epsilon,(n)}$ \textbf{does depend} on the data and samples. Algorithm \ref{alg:raoBlackwellizedFilter} lays out the steps for programming such a filter.
\begin{algorithm}
\caption{Rao-Blackwellized Particle Filter Algorithm at $(k+1)^{th}$ observation}
\label{alg:raoBlackwellizedFilter}
\begin{algorithmic}
\STATE $\Sigma_{k+1|k}=A\Sigma_kA^T+R^TR\qquad\qquad$ //compute the covariance of $\vec{\tilde X}_{k+1}^\epsilon$ given $\mathcal F_k$.
\STATE  $S_{k+1}=H\Sigma_{k+1|k}H^T+\widetilde\Delta t\qquad\qquad$ //compute covariance of  $Y_{k+1}^\epsilon-Y_k^\epsilon$ given $\mathcal F_k$
\STATE $K_k=\Sigma_{k+1|k}H^T/S_{k+1}\qquad\qquad$ //compute Kalman filter gain operator, 
\STATE $\Sigma_{k+1}=(I-K_kH)\Sigma_{k+1|k}\qquad\qquad$ //compute the covariance of $\vec{\tilde X}_{k+1}^\epsilon$ given $\mathcal F_{k+1}$.
\FOR{$n = 1 \to R$}
\FOR{$\ell=0\to m-1$}
\STATE $\widetilde\Theta_{mk+\ell+1}^{(n)}$ given $\widetilde\Theta_{mk+\ell}^{(n)}\qquad\qquad$ //sequentially update the particle
\ENDFOR
\STATE $\widehat{\vec{\widetilde X}}_{k+1|k}^{\epsilon,(n)}=A\widehat{\vec{\widetilde X}}_{k}^{\epsilon,(n)}+B\vec{\tilde\Theta}_{k+1}^{(n)}\qquad\qquad\qquad\qquad$ //predict the state at the next time step.
\STATE  $\omega_{k+1}^{(n)}=\exp\left\{-\frac{1}{2}\left(\frac{y_{k+1}-y_k-H\widehat{\vec{\widetilde X}}_{k+1|k}^{\epsilon,(n)}}{\sqrt{S_{k+1}}}\right)^2\right\}\cdot\omega_k^{(n)}\qquad\qquad$ //compute unnormalized weight
\STATE $\widehat{\vec{\widetilde X}}_{k+1}^{\epsilon,(n)}=\widehat{\vec{\widetilde X}}_{k+1|k}^{\epsilon,(n)}+K_k\left(y_{k+1}-y_k-H\widehat{\vec{\widetilde X}}_{k+1|k}^{\epsilon,(n)}\right)\qquad\qquad$ // update the filter.
\ENDFOR
\STATE $c_{k+1}=\sum_n\omega_{k+1}^{(n)}\qquad\qquad$ //compute normalization weight
\STATE $\omega_{k+1}=\omega_{k+1}/c_{k+1}\qquad\qquad$ //normalize
\STATE //////////////////////////////////////////////////////////////////////////////////////////
\STATE //        perform SIR if necessary
\IF{$1/\sum_n(\omega_{k+1}^{(n)})^2\leq \eta R$}
\STATE replace $\left\{\widetilde\Theta_{m(k+1)}^{(n)},\widehat{\vec{\widetilde X}}_{k+1}^{\epsilon,(n)}\right\}_{n=1}^R$ with an $\{\omega_{k+1}^{(n)}\}_{n=1}^R$-weighted bootstrap sample of size $R$
\STATE replace $\{\omega_{k+1}^{(n)}\}_{n=1}^R$ with $(1,1,\dots,1)/R$
\ENDIF
\end{algorithmic}
\end{algorithm}


\begin{thebibliography}{1}

\bibitem{asmussen2007} S. Asmussen, P. Glynn, \underline{Stochastic Simulation: Algorithms and Analysis}. Springer, 2007.

\bibitem{BCC97} B. Bakshi, C. Cao, Z. Chen (1997) ``Empirical Performance of Alternative Options Pricing Models," Journal of Finance, Vol. 52. No 5. pp. 2003-2049.


\bibitem{barShalom}Y. Bar-Shalom, X. R. LI, \underline{Estimation and Tracking: Principles, Techniques and
Software}. Boston: Artech House, 1993.

\bibitem{baum} L. Baum, T. Petrie, G. Soules, N. Weiss, ``A Maximization Technique Occurring in the Statistical Analysis of Probabilistic Functions of
Markov Chains," The Annals of Mathematical Statistics, Vol. 41, No. 1 (Feb., 1970), pp. 164-171.

\bibitem{baumWelch} L. Baum, L. Welch, ``Statistical Estimation Procedure for Probabilistic Functions of Finite Markov processes." Submitted for publication Proc. Nat. Acad. Sci. USA.

\bibitem{bhattBudhiraja} A. Bhatt, A. Buhiraja, A. Karandikar, ``Markov Property and Ergodicity of the Nonlinear Filter," SIAM J. CONTROL OPTIM. Vol. 39, No. 3, pp. 928Ð949
\bibitem{budhiraja} A. Budhiraja ``Ergodic Properties of Nonlinear Filter," Stochastic Processes and their Applications, Vol. 95, No. 1. (September 2001), pp. 1-24

\bibitem{cappe2005} O. Cappe, E. Moulines, T. Ryden,  \underline{Inference in Hidden Markov Models}. Springer 2005.


\bibitem{ethierKurtz} Stewart N. Ethier, Thomas G. Kurtz, \underline{Markov Processes: Characterization and Convergence}. Wiley, 1986.

\bibitem{FPS00} J.-P. Fouque, G. Papanicolaou, and R. Sircar ``Derivatives in Financial Markets with Stochastic Volatility", Cambridge University Press, 2000.

\bibitem{FPSS04} J.-P. Fouque, G. Papanicolaou, R. Sircar, K. S\o lna (2004) ``Maturity Cycles in Implied Volatility," Finance and Stochastics, vol. 8, pp. 451-477.

\bibitem{freidlin04} M. Freidlin, ``Some Remarks on the Kramers-Smoluchowski Approximation", Journal of Statistical Physics, 117, pages 617-634, 2004.


\bibitem{krylov}B.Fristedt, N. Jain, N. Krylov, \underline{Filtering and Prediction: A Primer}. American Mathematical Society, 2007.


\bibitem{gustafsson2002} F. Gustafsson, F. Gunnarsson,  N. Bergman, U. Forssell, J. Jansson, R. Karlsson, P. Nordlund, ``Particle Filters for Positioning, Navigation and Tracking,'' 2002.

\bibitem{howison} S. Howison, A. Rafailidis, H. Rasmussen, ``On the pricing and hedging of volatility derivatives''  Applied Mathematical Finance, Volume 11, Issue 4 December 2000 , pages 317 - 346.

\bibitem{jazwinski1970} A.H. Jazwinski, \underline{Stochastic Processes and Filtering Theory}. Academic Press, New York, 1970.


\bibitem{kushner08} H. Kushner ``Numerical Approximation to Optimal Nonlinear Filters," 2008, {\tt http://www.dam.brown.edu/lcds/publications/}

\bibitem{kutoyants} Y.A. Kutoyants. Statistical Inference for Ergodic Diffusion Processes. Springer, London, 2004

\bibitem{papanicolaou2010} A. Papanicolaou, ``Filtering for Fast Mean-Reverting Processes," Asymptotic Analysis, 70 (2010) 155-176.

\bibitem{pardoux} E. Pardoux, ``Filtrage Non Lineaire et Equations aux Derivees Partielles Stochastiques Associees," Ecol\'e d'\'et\'e de Probabilities de Saint-Flour, 1989.

\bibitem{perello} J. Perrel\'o, R. Sircar, J. Masoliver, ``Option Pricing Under Stochastic Volatility: the Exponential Ornstein-Uhlenbeck Model," Journal of Statistical Mechanics (2008), P06010.

\bibitem{rabiner1989} L. R. Rabiner ``A tutorial on Hidden Markov Models and selected applications in speech recognition". Proceedings of the IEEE 77 (2), (February 1989). pages: 257-286.

\bibitem{petrov1999} B. Rozovsky, A. Petrov, ``Optimal Nonlinear Filtering for Track-Before-Detect in IR Image Sequences.'' SPIE proceedings: Signal and Data Processing of Small Targets, vol. 3809, Denver, Co, 1999.

\bibitem{petrov2000} B. L. Rozovskii, A. Petrov, R. B. Blazek, ``Interacting Banks of Bayesian Matched Filters.'' SPIE Proceedings: Signal and Data Processing of Small Targets,
Vol. 4048, Orlando, FL, 2000. 

\bibitem{rozovsky1991} B. Rozovsky, ``A simple proof of uniqueness for Kushner and Zakai equations.''  In Stochastic analysis, ed. E. Mayer-Wolf, 449-458.  Boston: Academic Press, 1991. 

\bibitem{rozovskii1972} B. Rozovskii, A.N. Shiryaev, ``On Infinite Order Systems of Stochastic Differential Equations Arising in the Theory of Optimal Nonlinear Filtering," \textit{Theory of Probability and Applications} Volume 17, number 2, 1972.

\bibitem{gustafsson2005} T.  Sch\"{o}n, G. Gustafsson, P. Norlund, ``Marginalized Particle Filters for Mixed Linear/Non-Linear State-Space Models,''  Signal Processing, IEEE Transactions on
Volume 53, Issue 7, July 2005 Page(s): 2279 - 2289


\bibitem{yinZhangCont} G. George Yin, Qing Zhang, \underline{Continuous-Time Markov Chains and Applications: A Singular }\\
\underline{Perturbation Approach}. Springer, 1998.

\bibitem{yinZhangDiscrete} G. George Yin, Qing Zhang, \underline{Discrete-Time Markov Chains: Two-Time-Scale Methods and Applications}. Springer, 2005.

\bibitem{zakai1969} Zakai, M.  ``On the optimal filtering of diffusion processes.'' Zeit. Wahrsch. 11 230-243. 1969

\end{thebibliography}
\end{document}